\documentclass[11pt]{amsart}

\usepackage{epsfig, subfig, amsmath, amssymb, latexsym, amsfonts, xcolor, easybmat, bbm}

\newcommand\ol{\ensuremath{\overline}}

\newcommand\eop{{{\hfil \ensuremath \Box}}}

\newcommand\eps{\ensuremath {\varepsilon}}

\newenvironment{cor}{\subsection{}{\textbf {Corollary.}}\em}{}
\newenvironment{defn}{\subsection{}{\textbf {Definition.}}\em}{\smallskip}
\newenvironment{eg}{\subsection{}{\textbf {Example.}}}{\smallskip}

\newenvironment{lem}{\subsection{}{\textbf {Lemma.}}\em}{}
\newenvironment{notation}{\subsection{}{\textbf{Notation.}}}{\smallskip}

\newenvironment{prop}{\subsection{}{\textbf {Proposition.}}\em}{}

\newenvironment{rem}{\subsection{}{\textbf {Remark.}}}{\smallskip}
\newenvironment{thm}{\subsection{}{\textbf {Theorem.}}\em}{}

\newenvironment{pf}{\noindent{\textbf {Proof.}}} {\begin{flushright}\eop \end{flushright}\smallskip}

\newcommand\fC{\ensuremath{\mathfrak C}}

\newcommand\fE{\ensuremath{\mathfrak E}}

\newcommand\fP{\ensuremath{\mathfrak P}}

\newcommand\cB{\ensuremath{\mathcal B}}

\newcommand\cH{\ensuremath{\mathcal H}}

\newcommand\cK{\ensuremath{\mathcal K}}

\newcommand\cM{\ensuremath{\mathcal M}}
\newcommand\cN{\ensuremath{\mathcal N}}

\newcommand\cR{\ensuremath{\mathcal R}}

\newcommand\cU{\ensuremath{\mathcal U}}

\newcommand\bbC{\ensuremath{\mathbb C}}

\newcommand\bbM{\ensuremath{\mathbb M}}
\newcommand\bbN{\ensuremath{\mathbb N}}

\newcommand\bbR{\ensuremath{\mathbb R}}

\newcommand\bbZ{\ensuremath{\mathbb Z}}

\newcommand\ttt{\ensuremath{\textsc}}

\newcommand\bofh{\ensuremath{\cB ( \cH)}}
\newcommand\kofh{\ensuremath{\cK ( \cH)}}

\newcommand\hilb{\ensuremath{\mathcal H}}

\newcommand\norm{\ensuremath {\Vert}}

\newcommand\rank{\ensuremath{\mathrm{rank}\, }}

\newcounter{asst}
\newcounter{lab}
\newcounter{asstAA}
\newcounter{asstBB}
\newcounter{asstCC}
\newcounter{asstDD}
\newcounter{asstEE}
\newcounter{asstFF}
\newcounter{asstGG}
\newcounter{asstHH}
\newcounter{asstII}
\newcounter{asstJJ}
\newcounter{asstKK}
\newcounter{asstLL}
\newcounter{asstMM}
\newcounter{suppAA}
 \definecolor{myred}{rgb}{0.6,0,0}
 \definecolor{myblue}{rgb}{0,0.2,0.4}
\definecolor{mygreen}{rgb}{0.2,0.6, 0.5}
\begin{document}

%%%%%%%%%%%%%%%%%%%%%
%Title
%%%%%%%%%%%%%%%%%%%%%

\title{Around the closures of the set of commutators and  the set of differences of idempotent elements of $\bofh$}

%%%%%%%%%%%%%%%%%%%%%
% Authorship
%%%%%%%%%%%%%%%%%%%%%

%\thanks{${}^1$ Research supported in part by NSERC (Canada)}
%\thanks{${}^2$ Research supported in part by National Natural Science Foundation of China (No.: 12071174), Science and Technology
%Development Project of Jilin Province (No.: 2019013028JH)}
%

\thanks{{\ifcase\month\or Jan.\or Feb.\or March\or April\or May\or
June\or
July\or Aug.\or Sept.\or Oct.\or Nov.\or Dec.\fi\space \number\day,
\number\year}}
\author
	[L.W. Marcoux]{{Laurent W. Marcoux%${}^1$
}}
\address
	{Department of Pure Mathematics\\
	University of Waterloo\\
	Waterloo, Ontario \\
	Canada  \ \ \ N2L 3G1}
\email{Laurent.Marcoux@uwaterloo.ca}

%%%%%%%%

\author
	[H. Radjavi]{{Heydar Radjavi}}
\address
	{Department of Pure Mathematics\\
	University of Waterloo\\
	Waterloo, Ontario \\
	Canada  \ \ \ N2L 3G1}
\email{hradjavi@uwaterloo.ca}

%%%%%%%%

\author
	[Y.H.~Zhang]{{Yuanhang~Zhang%${}^2$
}}
\address
	{School of Mathematics\\
	Jilin University\\
	Changchun 130012\\
	P.R. CHINA}
\email{zhangyuanhang@jlu.edu.cn}

%
%\thanks{${}^1$ Research supported in part by NSERC (Canada)}
%\thanks{${}^2$ Research supported in part by National Natural Science Foundation of China (No.: 12071174)}

%%%%%%%%%%%%%%%%%%%%%
% Abstract
%%%%%%%%%%%%%%%%%%%%%

\begin{abstract}
We describe the norm-closures of  the set  $\fC_{\fE}$ of commutators of idempotent operators and the set $\fE - \fE$ of differences of idempotent operators acting on a finite-dimensional complex Hilbert space, as well as characterising the intersection of the closures of these sets with the set $\kofh$ of compact operators acting on an infinite-dimensional, separable Hilbert space.   Finally, we characterise the closures of the set $\fC_\fP$ of commutators of orthogonal projections and the set $\fP - \fP$ of differences of orthogonal projections acting on an arbitrary complex Hilbert space.
\end{abstract}

%%%%%%%%%%%%%%%%%%%%%
% Keywords and Subject Classification
%%%%%%%%%%%%%%%%%%%%%

\keywords{commutators, differences, idempotents, projections, closures}
\subjclass[2010]{Primary: 47B47. Secondary: 47A58}

\maketitle
\markboth{\textsc{  }}{\textsc{}}

%%%%%%%%%%%%%%%%%%%%%%%%%%%%%%%%%%%%%%%%%%
%%%%%%%%%%%%%%%%%%%%%%%%%%%%%%%%%%%%%%%%%%
% SECTION ONE
%%%%%%%%%%%%%%%%%%%%%%%%%%%%%%%%%%%%%%%%%%
%%%%%%%%%%%%%%%%%%%%%%%%%%%%%%%%%%%%%%%%%%

\section{Introduction}

%%%%%%%%%%%%%%%%%%%%%%%%%%%%%%%%%%%%%%%%%%

\subsection{} \label{sec1.01}
Let $\hilb$ be a complex Hilbert space.   By $\bofh$, we denote the norm-closed algebra of all bounded linear operators acting on $\hilb$.   There are surprisingly few well-understood classes of continuous linear operators acting on $\hilb$, but amongst the best understood of these are the orthogonal projections.   Recall that an operator $P \in \bofh$ is said to be an \textbf{orthogonal projection} if $P = P^* = P^2$.

The three most important notions of equivalence of Hilbert space operators are \textbf{similarity, unitary equivalence} and \textbf{approximate unitary equivalence}.    Given $A, B \in \bofh$, we shall write $A \sim B$ to indicate that $A$ is similar to $B$; i.e. there exists $S \in \bofh$ invertible such that $B = S^{-1} A S$.   We write $A \simeq B$ to indicate that $A$ is unitarily equivalent to $B$; i.e. there exists a unitary operator $U \in \bofh$ such that $B = U^* A U$.  Finally, we write $A \simeq_a B$ to indicate that $A$ is approximately unitarily equivalent to $B$, meaning that there exists a sequence $(U_n)_n$ of unitary operators in $\bofh$ such that $B = \lim_n U_n^* A U_n$.   This is equivalent to saying that the (norm) closure $\ttt{clos}\, \cU(A)$ of the unitary orbit $\cU(A) := \{ U^* A U: U \in \bofh \mbox{ unitary}\}$ of $A$ coincides with $\ttt{clos}\, \cU(B)$.

It is easy to verify that an operator $Q \in \bofh$ is approximately unitarily equivalent to a projection $P$ if and only if $Q$ is unitarily equivalent to $P$, in which case $Q$ is itself a projection.   Furthermore, it is a standard exercise in operator theory to show that an operator $E \in \bofh$ is similar to \emph{some} projection if and only if $E$ is \textbf{idempotent};  that is, $E^2 = E$.

There exists a substantial literature centred around the characterisation of specific linear and/or multiplicative combinations of projections and idempotents in $\bofh$, and indeed in other $C^*$-algebras~\cite{PearcyTopping1967, FongMurphy1985, Nishio1985,  Wu1990, GoldsteinPaskiewicz1992, Spiegel1993, LaurieMathesRadjavi1994, Wu1995, Shulman2003,  Rabanovych2004a, Rabanovych2004b,  Bikchentaev2005, KaftalNgZhang2012a, Andruchow2014, KomisarskiPaskiewicz2015}.

\smallskip

We shall focus on two particular instances of this problem, namely:  \emph{commutators} and \emph{differences}.    More specifically, our interest will lie in describing the \emph{norm-closures} of the sets described below.

%%%%%%%%%%%%%%%%%%%%%%%%%%%%%%%%%%%%%%%%%%

\subsection{Notation} \label{sec1.02} \ \ \
Given a Hilbert space $\hilb$, we define $\fE := \{ E \in \bofh: E = E^2 \}$, $\fP := \{ P \in \bofh: P = P^* = P^2\}$, and we set
\[
\fE - \fE := \{ E - F: E, F \in \fE\} \]
and
\[
\fP - \fP := \{ P - Q: P, Q \in \fP\}. \]
Following~\cite{HartwigPutcha1990}, we refer to differences of idempotents as $\ttt{doi}$s, and in light of this we refer to differences of projections as $\ttt{dop}$s.

We also define the sets
\begin{align*}
\fC_{\fE} & := \{ [E, F] : E, F \in  \fE \}  \mbox{\ \ \ and }\\
\mathfrak{C}_{\fP} &:= \{ [P, Q]:  P, Q  \in \fP\}.
\end{align*}
The elements of $\fC_\fE$ are commutators of idempotents (we shall refer to them as $\ttt{coi}$s), and we shall refer to elements of $\fC_\fP$ as $\ttt{cop}$s, short for \emph{commutators of projections}.

\smallskip

Finally (keeping in mind that an \textbf{involution} is an invertible operator $S \in \bofh$ such that  $S = S^{-1}$), we shall write
\begin{align*}
\ttt{Neg}_S (\hilb) &:= \{ T \in \bofh:  T \mbox{ is similar to } -T \}; \\
\ttt{Neg}_U (\hilb) &:= \{ T \in \bofh: T \mbox{ is unitarily equivalent to } -T \}; \mbox{ and} \\
\ttt{Neg}_{\ttt{invS}} (\hilb) &: = \{ T \in \bofh: T \mbox{ is involution-similar to } -T\}.
\end{align*}

Obviously
\[
\ttt{Neg}_U (\hilb) \subseteq \ttt{Neg}_S (\hilb) \mbox{\ \ \ \ and \ \ \ \ } \ttt{Neg}_{\ttt{inv}} (\hilb) \subseteq \ttt{Neg}_S (\hilb). \]

%%%%%%%%%%%%%%%%%%%%%%%%%%%%%%%%%%%%%%%%%%

\subsection{} \label{sec1.03} \ \ \ Our first goal will be to classify the norm-closures of the sets $\fC_\fE$ and $\fE - \fE$.   We note that in the case where the underlying Hilbert space is finite-dimensional, the (non-closed) sets themselves have been classified.

Indeed, when $\dim\, \hilb$ is finite, the characterisation of $\fC_\fE$  is due to Drnov\v{s}ek et al.~\cite[Theorem~8]{DrnovsekRadjaviRosenthal2002}.

%%%%%%%%%%%%%%%%%%%%%%%%%%%%%%%%%%%%%%%%%%

\begin{thm} \emph{\textbf{[Drnov\v{s}ek, Radjavi and Rosenthal.]}}\ \ \  \label{thm1.04}
If $n := \dim\, \hilb < \infty$ and $T \in \cB(\bbC^n)$, then $T \in \fC_\fE$ if and only if $T \sim -T$ and the Riesz component $T_1$ of $T$ corresponding to $\left\{ \frac{1}{2} i \right \}$ has the property that $T_1^2 + \frac{1}{4} I$ has a square root.
\end{thm}

%%%%%%%%%%%%%%%%%%%%%%%%%%%%%%%%%%%%%%%%%%

\smallskip

The corresponding result for $\fE - \fE$ is due to Hartwig and Putcha~\cite[Theorem~1b]{HartwigPutcha1990}.

\begin{thm} \emph{\textbf{[Hartwig and Putcha.]}}\ \ \  \label{thm1.05}
If $n := \dim\, \hilb < \infty$ and $T \in \cB(\bbC^n)$, then $T \in \fE-\fE$ if and only if the elementary divisors (see, e.g.~\cite{Jacobson1953}) of $T$ satisfy the following three conditions:
\begin{enumerate}
	\item[(i)]
	there are no restrictions on the elementary divisors $z^k$;
	\item[(ii)]
	the elementary divisors $(z-\alpha)^k$, $(z+\alpha)^k$ with $\alpha \ne 0, \pm 1$ occur in pairs with the same multiplicities; and
	\item[(iii)]
	the elementary divisors $(z-1)^{m_k}$, $(z+1)^{n_k}$, $k = 1, 2, \ldots, r$ obey $|m_k - n_k| \le 1$ when listed in non-increasing order.
\end{enumerate}	
\end{thm}

\bigskip

A complete characterisation of the sets $\fC_\fE$ and $\fE - \fE$ in the case where $\dim\, \hilb = \infty$ is not yet available, though the paper of Wang and Wu~\cite{WangWu1994} has many interesting partial results.  The problem of characterising $\ttt{clos}\, (\fC_\fE)$ and $\ttt{clos}\, (\fE - \fE)$ in this setting seem quite delicate.  To wit:  although all nilpotent operators of order two are known to lie in $\fC_\fE$, it is not known which nilpotent operators of order three are commutators of idempotents.

%%%%%%%%%%%%%%%%%%%%%%%%%%%%%%%%%%%%%%%%%%

\subsection{} \label{sec1.06}
Our study of the classes $\fC_\fE$ and $\fE - \fE$ will require us to to understand the spectrum of an operator which is similar to its own negative.   Recall that the \textbf{semi-Fredholm domain} of an operator $T \in \bofh$ is the set
\[
\rho_{sF}(T) := \{ \alpha \in \bbC: \ttt{ran}\, (T - \alpha I) \mbox{ is closed and } \min(\ttt{nul}(T - \alpha I), \ttt{nul} (T-\alpha I)^* ) < \infty \}. \]
A standard result (see, e.g.~\cite{CaradusPfaffenbergerYood1974}) shows that $\alpha \in \rho_{sF}(T)$ if and only if $\pi (T-\alpha I)$ is either left- or right-invertible in $\bofh/\kofh$, where $\kofh$ denotes the closed, two-sided ideal of compact operators acting on $\hilb$, and $\pi: \bofh \to \bofh/\kofh$ is the canonical quotient map.  When $\alpha \in \rho_{sF}(T)$, one defines the \textbf{semi-Fredholm index}
\[
\ttt{Ind}\, (T - \alpha I) = \ttt{nul}\, (T-\alpha I) - \ttt{nul}\, (T-\alpha I)^*. \]

Note that if $T \sim -T$, say $- T = R^{-1} T R$ for some invertible operator $R$, then $\alpha \in \rho_{sF}(T)$ implies and $\ttt{Ind} (T - \alpha I) = m \in \bbZ \cup \{ -\infty, \infty\}$ implies that
\[
R^{-1} (T - \alpha I) R = (-T - \alpha I) = - (T + \alpha I), \]
so that $T+\alpha I$ is semi-Fredholm and
\[
\ttt{Ind}\, (T+\alpha I) = \ttt{Ind}\, (- (T + \alpha I)) = \ttt{Ind}\, (T-\alpha I). \]

\smallskip

%%%%%%%%%%%%%%%%%%%%%%%%%%%%%%%%%%%%%%%%%%

In light of the above observations, the following definition will prove useful.

\smallskip

\begin{defn} \label{defn1.07}
Let $T\in \bofh$.  We say that $T$ is \textbf{balanced} if
\begin{enumerate}
	\item[(i)]
	$\sigma(T)=\sigma(-T)$;
	\item[(ii)]
	whenever  $\Omega_1, \Omega_2 \subseteq \bbC$ are disjoint open sets such that $\sigma(T) \subseteq \Omega_1 \cup 		\Omega_2$, then
		\[
		\dim \hilb({\Omega_1}; T)=\dim \hilb({-\Omega_1};T), \]
	where $\hilb({\Omega_1}; T)$ is the generalised eigenspace $($i.e. the range of the corresponding Riesz idempotent $E(\Omega_1; T))$ corresponding to
	$\sigma(T) \cap \Omega_1$;
	\item[(iii)]
%	if $\alpha$ is an eigenvalue of $T$ then so is $-\alpha$ and
%	\[
%	\ttt{nul}\, (T-\alpha I) = \ttt{nul}\, (T +\alpha I); \mbox{ and } \]
%	\item[(iv)]
	If $\alpha \in \bbC$, then $T - \alpha I$ is semi-Fredholm if and only if $T + \alpha I$ is semi-Fredholm, in which case
	\[
	\ttt{Ind}\, (T - \alpha I) = \ttt{Ind}\, (T + \alpha I). \]
\end{enumerate}	
We denote the set of balanced  operators by $\ttt{Bal}(\hilb)$.
\end{defn}

%%%%%%%%%%%%%%%%%%%%%%%%%%%%%%%%%%%%%%%%%%
%%%%%%%%%%%%%%%%%%%%%%%%%%%%%%%%%%%%%%%%%%
%%%%%%%%%%%%%%%%%%%%%%%%%%%%%%%%%%%%%%%%%%
%%%%%%%%%%%%%%%%%%%%%%%%%%%%%%%%%%%%%%%%%%
%%%%%%%%%%%%%%%%%%%%%%%%%%%%%%%%%%%%%%%%%%

\section{Elementary and general results} \label{sec2}

%%%%%%%%%%%%%%%%%%%%%%%%%%%%%%%%%%%%%%%%%%

\subsection{} \label{sec2.01}
Our ultimate goal would be to describe the relationships between each of the classes of operators defined above, as well as their norm-closures in $\bofh$.    In the case where $\dim\, \hilb = \infty$, our incomplete understanding of the sets $\fC_\fE$ and $\fE - \fE$ themselves complicates matters.  For this reason, when considering the closures of $\fC_\fE$ and $\fE - \fE$, in this paper we shall focus mostly on two cases.   First,  we shall direct our attention to the case where $\dim\, \hilb < \infty$, where the sets $\fC_\fE$ and $\fE - \fE$ are fully understood.   Next,  we turn to a description of the sets $\ttt{clos} (\fC_\fE) \cap \kofh$ and $\ttt{clos}\, (\fE - \fE) \cap \kofh$.

Because the sets $\fC_\fP$ and $\fP - \fP$ are fully understood independently of the dimension of the underlying Hilbert space, in Section~\ref{sec5} we shall be able to characterise their closures.

\smallskip

We begin with a couple of general and elementary observations concerning the classes $\fC_\fE$ and $\ttt{bal}(\hilb)$ which will be used throughout the paper.

%%%%%%%%%%%%%%%%%%%%%%%%%%%%%%%%%%%%%%%%%%

\begin{prop} \label{prop2.02}
For any Hilbert space $\hilb$,
\[
\fC_\fE \subseteq \ttt{Neg}_{\ttt{invS}} (\hilb) \subseteq \ttt{Neg}_S (\hilb) \subseteq \ttt{Bal} (\hilb). \]
\end{prop}

\begin{pf}
It is routine to verify that all four sets above are invariant under similarity, and this will be used implicitly below.

\smallskip

The proof that $\fC_\fE \subseteq \ttt{Neg}_{\ttt{invS}}$ is an immediate consequence of Theorem~1 in the paper~\cite{DrnovsekRadjaviRosenthal2002}, while the inclusion $\ttt{Neg}_{\ttt{invS}} (\hilb) \subseteq \ttt{Neg}_S (\hilb)$ is trivial.

\smallskip

If $T \in \ttt{Neg}_S (\hilb)$ and $-T = S^{-1} T S$, then clearly $\sigma(T) = \sigma(S^{-1} T S) = \sigma(-T)$,
\[
\dim\, \hilb(\Omega_1; T) = \dim\, \hilb(\Omega_1; S^{-1} T S) = \dim\, \hilb(\Omega_1; -T) = \dim\, \hilb(-\Omega_1; T), \]
and for $\alpha \in \rho_{sF}(T)$,
\[
\ttt{Ind}\, (T - \alpha I) = \ttt{Ind}\, S^{-1} (T-\alpha I) S = \ttt{Ind}\, (-T - \alpha I) = \ttt{Ind}\, (T+\alpha I), \]
so that $T \in \ttt{Bal} (\hilb)$.

\end{pf}

%%%%%%%%%%%%%%%%%%%%%%%%%%%%%%%%%%%%%%%%%%

The following observation can easily be derived from the above result and Lemma~1.1 of the paper of Wang and Wu~\cite{WangWu1994}.  We include an alternate proof since it is so short.

\begin{prop} \label{prop2.03}
For any Hilbert space $\hilb$,
\[
\fC_\fE \subseteq \fE - \fE. \]
\end{prop}

\begin{pf}
Since both of these sets are invariant under similarity, it suffices to consider $T \in \fC_\fE$ of the form $T = [P, F]$, where $P = \begin{bmatrix} I & 0 \\ 0 & 0 \end{bmatrix}$ is an orthogonal projection, and $F  = \begin{bmatrix} F_1 & F_2 \\ F_3 & F_4 \end{bmatrix}$ relative to this decomposition.

Then
\[
T = [P, F] = \begin{bmatrix} 0 & F_2 \\ -F_3 & 0 \end{bmatrix} = \begin{bmatrix} I & F_2 \\ 0 & 0 \end{bmatrix} - \begin{bmatrix} I & 0 \\ F_3 & 0 \end{bmatrix} \in \fE - \fE. \]

\end{pf}

%%%%%%%%%%%%%%%%%%%%%%%%%%%%%%%%%%%%%%%%%%

\subsection{} \label{sec2.04}  We finish this section with some easy observations.
\begin{itemize}
	\item{}
	$\fC_\fE$ is self-adjoint:   note that $[E, F]^* = (E F - F E)^* = F^* E^* - E^* F^* = [F^*, E^*]$ and $E^*, F^*$ are idempotents when $E, F$ are.
	\item{}
	From this it follows that $\ttt{clos}\, (\fC_\fE)$ is also self-adjoint.
	\item{}
	If $T \in \fC_\fP$, then $H := i T$ is self-adjoint, so that $\ttt{clos}\,  (i \fC_\fP)$ is also contained in the set $\bofh_{sa}$ of self-adjoint operators on $\hilb$.
	\item{}
	The set $\fP - \fP$ and its closure $\ttt{clos}\, (\fP - \fP)$ are contained in $\bofh_{sa}$.
\end{itemize}

%%%%%%%%%%%%%%%%%%%%%%%%%%%%%%%%%%%%%%%%%%
%%%%%%%%%%%%%%%%%%%%%%%%%%%%%%%%%%%%%%%%%%
%%%%%%%%%%%%%%%%%%%%%%%%%%%%%%%%%%%%%%%%%%
%%%%%%%%%%%%%%%%%%%%%%%%%%%%%%%%%%%%%%%%%%
%%%%%%%%%%%%%%%%%%%%%%%%%%%%%%%%%%%%%%%%%%
%%%%%%%%%%%%%%%%%%%%%%%%%%%%%%%%%%%%%%%%%%

\section{The closures of $\fC_\fE$ and of $\fE - \fE$ in the finite-dimensional setting}

\bigskip

%%%%%%%%%%%%%%%%%%%%%%%%%%%%%%%%%%%%%%%%%%

\subsection{} \label{sec3.01}
We now turn our attention to the case where $n := \dim\, \hilb < \infty$, and concentrate on the problem of describing the closures of the set $\fC_\fE$ of commutators of idempotent operators and the set $\fE - \fE$ of differences of idempotent operators in $\cB(\bbC^n) \simeq \bbM_n(\bbC)$.
	
%%%%%%%%%%%%%%%%%%%%%%%%%%%%%%%%%%%%%%%%%%

\begin{prop}  \label{prop3.02}
If $\dim\, \hilb < \infty$, then $\ttt{Neg}_U (\hilb)$ is closed.
\end{prop}

\begin{pf}

Suppose that $(T_m)_m$ is a sequence in $\ttt{Neg}_U(\hilb)$, and that $T = \lim_m T_m$.

\smallskip

Choose $U_{m}$ unitary such that $-T_m = U_m^* T_m U_m$.   Since $(U_m)_m$ is bounded, there exists a subsequence $(U_{m_j})_j$ which converges in norm to a (necessarily unitary operator) $V \in \bofh$.

Then
\[
V^* T V = \lim_{j} U_{m_j}^* T_{m_j} U_{m_j} = \lim_j - T_{m_j} = - T, \]
so that $T \in \ttt{Neg}_U(\hilb)$.
\end{pf}
	
%%%%%%%%%%%%%%%%%%%%%%%%%%%%%%%%%%%%%%%%%%

\begin{notation} \label{not3.03}
If $\dim\, \hilb < \infty$, then given an operator $T \in \bofh$ and an eigenvalue $\alpha \in \sigma(T)$, we denote by  $\mu(\alpha)$ the algebraic multiplicity of $\alpha$.
\end{notation}

%%%%%%%%%%%%%%%%%%%%%%%%%%%%%%%%%%%%%%%%%%

\begin{prop} \label{prop3.04}
If $\dim\, \hilb = n < \infty$, then
\[
\ttt{Bal} (\hilb) = \{ T \in \bofh: \alpha \in \sigma(T) \mbox{ implies } - \alpha \in \sigma(T) \mbox{ and } \mu(\alpha) = \mu(-\alpha) \}, \]
and thus $\ttt{Bal} (\hilb)$ is closed.
\end{prop}

\begin{pf}
Let $T \in \ttt{Bal}(\hilb)$.  Since $\sigma(T) = \sigma(-T)$ by definition of $\ttt{Bal} (\hilb)$, it follows that $\alpha \in \sigma(T)$ implies that $-\alpha \in \sigma(T)$.   Also, taking $\Omega_\eps = \{ z \in \bbC: | z - \alpha| < \eps \}$ for sufficiently small $\eps > 0$ (to ensure that $\Omega_\eps \cap \sigma(T) = \{ \alpha\}$), condition (ii) from Definition~\ref{defn1.07} implies that $\mu(\alpha) = \mu(-\alpha)$.   Thus
\[
\ttt{Bal} (\hilb) \subseteq \{ T \in \bofh: \alpha \in \sigma(T) \mbox{ implies } - \alpha \in \sigma(T) \mbox{ and } \mu(\alpha) = \mu(-\alpha) \}. \]

\

Conversely, the condition that $\alpha \in \sigma(T)$ implies $-\alpha \in \sigma(T)$ is equivalent to the statement that $\sigma(T) = \sigma(-T)$, and if $\Omega \subseteq \bbC$ is any open set which non-trivially intersects $\sigma(T)$, then there exist  $1 \le \kappa \le n$ and $\alpha_1, \alpha_2, \ldots, \alpha_\kappa \in \sigma(T)$ such that $\Omega \cap \sigma(T) = \{ \alpha_1, \alpha_2, \ldots, \alpha_\kappa\}$.   Thus
\[
\dim\, \hilb(\Omega; T) = \sum_{j=1}^\kappa \mu(\alpha_j) = \sum_{j=1}^\kappa \mu(-\alpha_j) = \dim\, \hilb(-\Omega; T). \]

\

Note that condition (iii) from Definition~\ref{defn1.07} always holds in the finite-dimensional setting.   Thus
\[
\{ T \in \bofh: \alpha \in \sigma(T) \mbox{ implies } - \alpha \in \sigma(T) \mbox{ and } \mu(\alpha) = \mu(-\alpha) \} \subseteq \ttt{Bal} (\hilb), \]
so that equality of these two sets holds.

Suppose that $(T_m)_m$ is a sequence in
\[
\{ T \in \bofh: \alpha \in \sigma(T) \mbox{ implies } - \alpha \in \sigma(T) \mbox{ and } \mu(\alpha) = \mu(-\alpha) \}\]
and that $T = \lim_m T_m$.
Since the function $\sigma$ that takes an element $T \in \bofh$ to its spectrum $\sigma(T) \subseteq \bbC$ is continuous when $\dim\, \hilb$ is finite, we see that $\alpha \in \sigma(T)$ implies that $-\alpha \in \sigma(T)$.   Furthermore, if $\Omega$ is any open neighbourhood of $\alpha \in \sigma(T)$ such that $\Omega \cap \sigma(T) = \{ \alpha\}$, then for all $m \ge 1$,
\[
\dim\, \hilb(\Omega; T_m) = \dim\, \hilb(-\Omega; T_m) \]
since $T_m$ is balanced, and so
\[
\dim\, \hilb(\Omega; T) = \dim\, \hilb(-\Omega; T). \]
From this it follows that
\[
\{ T \in \bofh: \alpha \in \sigma(T) \mbox{ implies } - \alpha \in \sigma(T) \mbox{ and } \mu(\alpha) = \mu(-\alpha) \}\]
is closed, and thus that $\ttt{Bal} (\hilb)$ is closed.
\end{pf}

%%%%%%%%%%%%%%%%%%%%%%%%%%%%%%%%%%%%%%%%%%

\smallskip

It is worth observing that as a consequence of Proposition~\ref{prop3.04}, an operator $T \in \cB(\bbC^n)$ is balanced if and only if its characteristic polynomial is either an even function or an odd function.

\bigskip

%%%%%%%%%%%%%%%%%%%%%%%%%%%%%%%%%%%%%%%%%%

\begin{thm} \label{thm3.05}
Suppose that $\dim\, \hilb < \infty$.   Then
\[
\ttt{clos}\, (\fC_\fE) = \ttt{Bal}(\hilb). \]
\end{thm}

\begin{pf}
	By Proposition~\ref{prop2.02}, $\fC_\fE \subseteq \ttt{Bal} (\hilb)$.  Thus
	\[
	\ttt{clos}\, (\fC_\fE) \subseteq \ttt{clos}\, (\ttt{Bal} (\hilb)) = \ttt{Bal} (\hilb). \]

\smallskip

	By Proposition~\ref{prop3.04}, we have that
	\[
	\ttt{Bal} (\hilb) = \{ T \in \bofh: \alpha \in \sigma(T) \mbox{ implies } - \alpha \in \sigma(T) \mbox{ and } \mu(\alpha) = \mu(-\alpha) \}. \]

	Now consider the converse, and suppose that $T \in \ttt{Bal}(\hilb)$.  Then we may write the elements of $\sigma(T)$ (repeated according to their algebraic multiplicity) as an $n$-tuple
	\[
	\Sigma_T := (\alpha_1, \alpha_2, \ldots, \alpha_n), \]
	where $\alpha_{2 k} = - \alpha_{2k-1}$ for all $1 \le k \le \frac{n}{2}$, and $\alpha_n = 0$ if $n$ is odd.  Furthermore, we can upper-triangularise $T$ so that if $[T] = [t_{i j}]$, then
\begin{itemize}
	\item{}
	$t_{i j} = 0$ if $i > j$, and
	\item{}
	$t_{k k} = \alpha_k$ for all $1 \le k \le n$.
\end{itemize}	

\smallskip

Let $\eps > 0$.   It is relatively easy to show that we can find $\beta_k \in \bbC$, $1 \le k \le n$ such that
\begin{itemize}
	\item{}
	each $| \beta_k - \alpha_k| < \eps$;
	\item{}
	if $i \ne j$, then $\beta_i \ne \beta_j$;
	\item{}
	for all $1 \le k \le \frac{n}{2}$, $\beta_{2k} = -\beta_{2k -1}$; and
	\item{}
	$\beta_n = 0$ if $n$ is odd.
\end{itemize}

Let $D = \ttt{Diag} (\beta_1, \beta_2, \ldots, \beta_n)$.   Then $D$ is a normal operator and (by the last two conditions) $D$ is unitarily equivalent to $-D$.
Furthermore, all eigenvalues of $D$ are distinct.    By considering Jordan forms, any operator $X \in \bbM_n(\bbC)$ such that $\sigma(X) = \sigma(D)$ is similar to $D$.

Let $X \in \bbM_n(\bbC)$ be the operator whose matrix is $[x_{i j}]$, where
\[
x_{i j} = \begin{cases} t_{i j} & \mbox{ if } i \ne j \\ \beta_k & \mbox{ if } i = k = j. \end{cases} \]
Then $X$ is similar to $D$.

Now, by Proposition~3 of~\cite{DrnovsekRadjaviRosenthal2002}, $D \in \fC_{\fE}$.   Since $\fC_{\fE}$ is invariant under similarity, $X \in \fC_{\fE}$.    But
\begin{align*}
	X - T
		&= \ttt{Diag} (x_{11} - t_{11}, x_{22} - t_{22}, \ldots, x_{nn} - t_{nn}) \\
		&= \ttt{Diag} (\beta_1-\alpha_1, \beta_2 - \alpha_2, \ldots, \beta_n - \alpha_n),
\end{align*}		
so $\norm X - T \norm = \max_{1 \le k \le n} |\beta_k - \alpha_k| < \eps.$

Since $\eps > 0$ is arbitrary, it follows that $T \in \ttt{clos}\, (\fC_{\fE})$.
\end{pf}
	
%%%%%%%%%%%%%%%%%%%%%%%%%%%%%%%%%%%%%%%%%%

\begin{cor} \label{cor3.06}
Suppose that $\dim\, \hilb < \infty$.   Then
\[
\ttt{clos}\, (\fC_\fE) = \ttt{clos}\, (\ttt{Neg}_{\ttt{invS}} (\hilb)) =\ttt{clos}\,  (\ttt{Neg}_S (\hilb)) = \ttt{clos}\, (\ttt{Bal} (\hilb)) = \ttt{Bal} (\hilb). \]
\end{cor}

\begin{pf}
This is an immediate consequence of Proposition~\ref{prop2.02}, Proposition~\ref{prop3.04} and Theorem~\ref{thm3.05}.
\end{pf}

%%%%%%%%%%%%%%%%%%%%%%%%%%%%%%%%%%%%%%%%%%

We now turn our attention to the closure of the set $\fE - \fE$ of differences of idempotents.

%%%%%%%%%%%%%%%%%%%%%%%%%%%%%%%%%%%%%%%%%%

\begin{lem} \label{lem3.07}
\begin{enumerate}
	\item[(a)]
	Let $\alpha \in \bbC$.   Then $B := \alpha \oplus -\alpha \in \cB(\bbC^2)$ is a difference of idempotents.
	\item[(b)]
	If $\dim\, \hilb < \infty$ and $N \in \bofh$ is  nilpotent, then $N \in \fE - \fE$.
	\item[(c)]
	If $\hilb_k$ is a Hilbert space and $B_k \in \cB(\hilb_k)$ lies in $\fE - \fE$, $1 \le k \le K$, then $B := \oplus_{k=1}^K B_k$ is a difference of idempotents in $\cB(\oplus_{k=1}^K \hilb_k)$.
\end{enumerate}
\end{lem}
	
\begin{pf}
\begin{enumerate}
	\item[(a)]
	If $\alpha = 0$, then $B = 0$ is trivially a difference of idempotents.   If $\alpha \ne 0$, then
	\[
	B \simeq \begin{bmatrix} 0 & \alpha \\ \alpha & 0 \end{bmatrix} = \begin{bmatrix} 1 & \alpha \\ 0 & 0 \end{bmatrix} - \begin{bmatrix} 1 & 0 \\ - \alpha & 0 \end{bmatrix} \in \fE - \fE. \]
	\item[(b)]
	By Proposition~6 of \cite{DrnovsekRadjaviRosenthal2002}, $N \in \fC_\fE \subseteq \fE - \fE$.
	\item[(c)]
	Setting $B_k = [E_k, F_k]$, $1 \le k \le K$ and $E := \oplus_{k=1}^K E_k$, $F := \oplus_{k=1}^K F_k$, we see that $B = [E, F]$.
\end{enumerate}
\end{pf}
	
%%%%%%%%%%%%%%%%%%%%%%%%%%%%%%%%%%%%%%%%%%

\begin{prop} \label{prop3.08}
Suppose that $2 \le n:= \dim \, \hilb < \infty$.   If $T \in \ttt{Bal}(\hilb)$, then $T \in \ttt{clos} \,(\fE - \fE)$.   Nevertheless, there exists $T \in \ttt{Bal}(\hilb)$ such that $T \not \in \fE - \fE$.
\end{prop}

\begin{pf}
By Theorem~\ref{thm3.05}, $\ttt{Bal}(\hilb) = \ttt{clos}\, (\fC_\fE)$, while Proposition~\ref{prop2.03} shows that $\fC_\fE \subseteq \fE - \fE$, whence
\[
\ttt{Bal}(\hilb) \subseteq \ttt{clos}\, (\fE - \fE). \]

\smallskip

As for the second statement, note that if $X = \begin{bmatrix} 2 & 0 \\ 0 & 2 \end{bmatrix}$ and $Y = \begin{bmatrix} - 2 & 1 \\ 0 & -2 \end{bmatrix}$, then $W = X \oplus Y$ is invertible, as are $W + I_4$ and $W-I_4$.    Note that $W \in \ttt{Bal}(\bbC^4)$.   By Proposition~4 of \cite{HartwigPutcha1990}, $W \in \fE - \fE$ if and only if $W$ is similar to $D \oplus - D$ for some invertible operator $D \in \cB(\bbC^2)$.   But clearly the Jordan form of $W$ prohibits this from happening.   Thus $W \in \ttt{Bal}(\bbC^4) \setminus (\fE - \fE)$.
\end{pf}

%%%%%%%%%%%%%%%%%%%%%%%%%%%%%%%%%%%%%%%%%%

\begin{prop} \label{prop3.09}
Let $2 \le n := \dim\, \hilb < \infty$, and let $\mathfrak{sl}(\hilb) := \{ A \in \bofh: \ttt{tr}(A) = 0\}$.    Then
\[
\mathfrak{sl}(\hilb) \cap \ttt{clos} (\fE-\fE) = \ttt{Bal}(\hilb). \]
\end{prop}

\begin{pf}
Clearly $\ttt{Bal}(\hilb) \subseteq \mathfrak{sl}(\hilb)$, and by Proposition~\ref{prop3.08}, we have that $\ttt{Bal}(\hilb) \subseteq \ttt{clos} (\fE - \fE)$, so that
\[
\ttt{Bal}(\hilb) \subseteq \mathfrak{sl}(\hilb) \cap \ttt{clos} (\fE-\fE). \]
Now suppose that $T \in \mathfrak{sl}(\hilb) \cap (\fE-\fE)$.   We may decompose $\hilb$ as $\hilb = \hilb_1 \oplus \hilb_2 \oplus \hilb_3$ in such a way that relative to this decomposition, we have
\[
T = \begin{bmatrix} T_{1 1} & T_{1 2} & T_{1 3} \\ 0 & T_{2 2} & T_{2 3} \\ 0 & 0 & T_{3 3} \end{bmatrix}, \]
where $\sigma(T_{1 1}) = \{ 0\}$, $\sigma(T_{2 2}) \subseteq \{ -1, 1\}$ and $\sigma(T_{3 3}) \subseteq \bbC \setminus \{ -1, 0, 1\}$.    Since the spectra of $T_{1 1}, T_{2 2}$ and $T_{3 3}$ are all disjoint, we see that $T$ is similar to $T_{1 1} \oplus T_{2 2} \oplus T_{3 3}$, and that we find from Lemma~2 of \cite{HartwigPutcha1990} that $T \in \fE - \fE$ if and only if each of $T_{1 1}, T_{2 2}$ and $T_{3 3}$ is.

Now $\sigma(T_{1 1}) = \{ 0\}$, so $T_{1 1}$ is balanced.   Also,  by Proposition~4 of \cite{HartwigPutcha1990}, $T_{3 3} \in \fE - \fE$ if and only if $T_{3 3}$ is similar to $D \oplus - D$ for some invertible operator $D \in \cB(\hilb_3)$, implying that $T_{3 3}$ is balanced.     Thus $\ttt{tr}(T_{1 1}) = 0 = \ttt{tr}(T_{3 3})$.   Since $\ttt{tr}(T) = 0$, it follows that $\ttt{tr}(T_{2 2}) = 0$.   But $\sigma(T_{2 2}) \subseteq \{ -1, 1\}$, which then implies that $T_{2 2}$ is balanced.     Hence $T \in \ttt{Bal}(\hilb)$, being the direct sum of balanced operators.

Finally, if $X \in \mathfrak{sl}(\hilb) \cap \ttt{clos}(\fE-\fE)$, then $X = \lim_n T_n$, where each $T_n \in \fE - \fE$.   Since $\ttt{tr}(\cdot)$ is continuous and integer-valued on $\fE - \fE$, and since $\ttt{tr}(X) = 0$, it follows that there exists $n_0 \ge 1$ such that $T_n \in \mathfrak{sl}(\hilb) \cap (\fE- \fE)$ for all $n \ge n_0$.    But then $T_n \in \ttt{Bal}(\hilb)$ for all $n \ge n_0$.   Since $\ttt{Bal}(\hilb)$ is closed by Proposition~\ref{prop3.04}, $X \in \ttt{Bal}(\hilb)$.  This completes the proof.
\end{pf}

%%%%%%%%%%%%%%%%%%%%%%%%%%%%%%%%%%%%%%%%%%

Let $n \ge 2$ and denote by $\{ e_k: 1 \le k \le n\}$ the standard orthonormal basis for $\bbC^n$.   We shall denote by $J_n$ the standard $n \times n$ Jordan cell in $\bbM_n(\bbC)$;  that is, $J_n$ is the unique operator satisfying  $J_n e_1 = 0$, while $J_n e_k = e_{k-1}$, $2 \le k \le n$.

\bigskip

%%%%%%%%%%%%%%%%%%%%%%%%%%%%%%%%%%%%%%%%%%

\begin{lem}\label{lem3.10}
Let  $\hilb$ be a Hilbert space with $n = \dim\, \hilb < \infty$  and $Z \in \bofh \simeq \bbM_n(\bbC)$.  Suppose that
\begin{enumerate}
	\item[(i)]
	$\sigma(Z) \subseteq \{ -1, 1\}$,
	\item[(ii)]
	$\ttt{tr}(Z) = r \in \bbN$; and
	\item[(iii)]
	$\ttt{nul}\, (Z-I) \ge r$.
\end{enumerate}
Then $Z \in \ttt{clos} (\fE - \fE)$.
\end{lem}

\begin{pf}
Suppose first that  $\sigma(Z)=\{1\}$. Then, from (ii) and (iii), we know
each elementary divisor of $Z$ has degree one. Hence, $Z$ is similar to $I_r$, whence $Z=I_r$.
Thus we may assume that $\sigma(Z)=\{-1,1\}$.

Consider the Jordan form of $Z$, namely
\[
	Z \simeq \left [ \oplus_{j=1}^{\kappa^+} (I_{m_j} + J_{m_j}) \right] \oplus \left[ \oplus_{j=1}^{\kappa^-} (-I_{n_j} + J_{n_j}) \right]. \]
Let $s := \sum_{j=1}^{\kappa^-} n_j$. Observe that
\begin{itemize}
	\item{}
	the fact that $\ttt{nul}\, (Z-I) \ge r$ implies that $\kappa^+ \ge r$; and
	\item{}
	the fact that $\ttt{tr}(Z) = r$ implies that $\sum_{j=1}^{\kappa^+} m_j - s = r$.
%    the fact that $Z$ has no degree one elementary divisors implies that  $s\geq r$.
\end{itemize}
Here we agree that $J_1=0$.	

Since $-1\in \sigma(Z)$, $s\geq 1$. If $s<r$, set $c=r-s$. According to the above facts, it is easy to deduce that there exists a subset $\Lambda\subset \{1,\cdots,\kappa^+\}$ with $|\Lambda|=c$,
such that $m_j=1$, when $j\in \Lambda$. By reindexing, we may assume that $\Lambda=\{\kappa^+-c+1,\cdots, \kappa^+\}$. Then
\[Z \simeq \left [ \oplus_{j=1}^{\kappa^+-c} (I_{m_j} + J_{m_j}) \right] \oplus \left[ \oplus_{j=1}^{\kappa^-} (-I_{n_j} + J_{n_j}) \right]\oplus I_c\doteq  Z_1\oplus I_c.\]
It is clear that $Z \in \ttt{clos} (\fE - \fE)$ if $Z_1\in \ttt{clos} (\fE - \fE)$. Now let $r_1=s$, then $\sigma(Z_1)=\{-1,1\}$, $\ttt{tr}(Z_1) = r_1 \in \bbN$, and $\ttt{nul}\, (Z_1-I) \ge r_1$.
Furthermore, $s_1:= \sum_{j=1}^{\kappa^-} n_j=s=r_1$. Hence, without loss of generality, we may add a further assumption that ``$s\geq r$" to $Z$.

For each $n \ge 1$, we define a diagonal operator
\[
	D_n = \ttt{diag} (d_1^{(n)}, d_2^{(n)}, \ldots, d_s^{(n)}) \]
with the properties that
\begin{itemize}
	\item[(i)]
	for any fixed $n \ge 1$, all of the diagonal entries $d_j^{(n)} \in (0, 1)$, $1 \le j \le s$ are distinct, and
	\item[(ii)]
	$\lim_n d_j^{(n)} = 1$ for all $1 \le j \le s$.
\end{itemize}
Partition the set $\{ 1, 2, \ldots, s \}$ into $r+1$ disjoint sets $\Omega_1, \Omega_2, \ldots, \Omega_{r+1}$, where $| \Omega_{j}| = m_j-1$ for $1 \le j \le r$ and $\Omega_{r+1}$ contains the remaining elements of $\{ 1, 2, \ldots, s\}$.  (It is possible that $\Omega_{r+1}$ might be empty.)    Define $D_j^{(n)} = \ttt{diag} \{ d_\ell^{(n)} : \ell \in \Omega_j\}$, $1 \le j \le r+1$, and note that $D_n \simeq \oplus_{j=1}^{r+1} D_j^{(n)}$.   Set
\[
Y_n^+ = \left[ \oplus_{j=1}^{r} (1 \oplus D_j^{(n)}) \right] \oplus  D_{r+1}^{(n)}. \]
(The point is that $Y_n^+$ is a direct sum of $r$ diagonal operators acting on spaces of dimension $m_1, m_2, \ldots, m_{r}$, and each of these diagonal operators has first entry equal to $1$, along with another diagonal operator which brings the dimension of the space upon which $Y_n^+$ acts up to $s$.  Other than the $1$'s which appear $r := \ttt{tr}(Z)$ times, all other diagonal entries of $Y_n^+$ should be distinct.)
Because all of the diagonal entries of $(1 \oplus D_j^{(n)})$ are distinct, we see that
\[
(1 \oplus D_j^{(n)}) \sim A_j^{(n)} := (1 \oplus D_j^{(n)}) + J_{m_j}, \ \ \  1 \le j \le r, \ \ n \ge 1. \]
Also, because all of the diagonal entries of $D_{r+1}^{(n)}$ are distinct (for any $n \ge 1$), we see that
\[
D_{r+1}^{(n)} \sim A_{r+1}^{(n)} := D_{r+1}^{(n)} + \left[ \oplus_{j=r+1}^{\kappa^+} J_{m_j} \right]. \]
In other words,
\[
Y_n^+ \sim B_n := \left[ \oplus_{j=1}^r A_j^{(n)} \right] \oplus A_{r+1}^{(n)}, \ \ \ n \ge 1. \]
Observe also that
\[
\lim_n B_n = \lim_n \left[ \oplus_{j=1}^r A_j^{(n)} \right] \oplus A_{r+1}^{(n)} = \left[ \oplus_{j=1}^{\kappa^+} (I_{m_j} + J_{m_j}) \right]. \]

\bigskip

Next, let $Y_n^- := - D_{n}$.   Since all of the diagonal entries of $D_n$ were distinct (for all $n \ge 1$), we see that
\[
Y_n^- = -D_n \sim C_n :=  -D_n + \left[ \oplus_{j=1}^{\kappa^-} J_{n_j} \right], \]
and thus
\[
\lim_n C_n = \left[ \oplus_{j=1}^{\kappa^-} (-I_{n_j} + J_{n_j}) \right]. \]

But $Y_n := Y_n^+ \oplus Y_n^- \simeq I_r \oplus D_n \oplus -D_n \in \fE-\fE$ by the Hartwig-Putcha Theorem~\cite[Theorem 1a]{HartwigPutcha1990}, and $Y_n \sim B_n \oplus C_n$, implying that $B_n \oplus C_n \in \fE-\fE$.  Finally,
\[
\lim_n (B_n \oplus C_n) = Z, \]
and therefore $Z \in \ttt{clos} (\fE-\fE)$.
\end{pf}

%%%%%%%%%%%%%%%%%%%%%%%%%%%%%%%%%%%%%%%%%%

\begin{lem}\label{lem3.11}
Let  $\hilb$ be a Hilbert space with $n = \dim\, \hilb < \infty$  and let $T \in \cB(\bbC^n) \simeq \bbM_n(\bbC)$.   Suppose that $\ttt{tr}(T) = r \in \bbN$.
If $T \in \ttt{clos} (\fE - \fE)$, then $\ttt{nul}\, (T-I) \ge r$.
\end{lem}

\begin{pf}
Suppose that $T=\underset{n\to \infty}{\lim} T_n$, where $T_n\in \fE-\fE$.  By the continuity of the spectrum and therefore of the
trace, we may assume without loss of generality that $\ttt{tr}(T_n) = r$ for all $n \in \mathbb{N}$.

By the Hartwig-Putcha Theorem~\ref{thm1.05}, the eigenvalues of $T_n$ which belong to $\mathbb{C}\setminus \{-1,1\}$ come in pairs $\{-\alpha,\alpha\}$, and therefore do
not contribute to the trace. From this, and again by the Hartwig-Putcha Theorem, it follows that if we fix the exponents $m_j^{(n)}$
occurring
in the elementary divisors of $T_n$ corresponding to 1 and $n_j^{(n)}$ corresponding to $-1$, for any
$n\in \mathbb{N}$, we will always have at least $r$ $j$'s for which
\[m_j^{(n)}-n_j^{(n)}=1.\]
It turns out that $\ttt{nul}\, (T_n-I) \ge r$.
Since the function $\ttt{nul}(\cdot)$ is  upper-semicontinuous (for the rank function is lower-semicontinuous),
$\ttt{nul}\, (T-I) \ge r$.
\end{pf}

%%%%%%%%%%%%%%%%%%%%%%%%%%%%%%%%%%%%%%%%%%

The hypothesis in the next theorem that the trace of the operator $T$ should be non-negative is there only to simplify the statement of the result.   Note that $T \in \ttt{clos} (\fE - \fE)$ if and only if $-T \in \ttt{clos} (\fE-\fE)$, so by replacing $T$ by $-T$ if necessary, the trace of $T$ may always be assumed to be a non-negative integer.

\smallskip
	
\begin{thm} \label{thm3.12}
Let $T \in \cB(\bbC^n) \simeq \bbM_n(\bbC)$ and suppose that $\ttt{tr} (T) = r \in \bbN$.  The following are equivalent.
\begin{enumerate}
	\item[(a)]
	$T \in \ttt{clos}(\fE - \fE)$; and
	\item[(b)]
	$T \sim B \oplus Z$, where $B$ is balanced,  $\sigma(Z) \subseteq \{ -1, 1\}$ and $\ttt{nul}\, (Z-I) \ge r$.
\end{enumerate}	
\end{thm}

\begin{pf}
\begin{enumerate}
	\item[(a)] implies (b). \ \ \ By using Riesz idempotent theorem, we may assume that $T=B\oplus Z$, where $\sigma(B)\in \mathbb{C}\setminus\{-1,1\}$, and $\sigma(Z)\subset\{-1,1\}$.
	We claim that $B$ is balanced.  Indeed, there exists $0 < \delta < 1$ such that $\sigma(B) \subseteq \{ z \in \bbC: |z| < \delta\}$.    If $T = \lim_m T_m$ where $T_m \in \fE - \fE$ for all $m \ge 1$, then by the continuity of the map $X \mapsto \sigma(X)$ in the finite-dimensional setting (with eigenvalues counted according to their algebraic multiplicities), we see that $\sigma(T)$ is the limit of $\sigma(T_m)$, and thus $\sigma(B)$ is the limit of $\sigma(T_m) \cap \{ z \in \bbC: |z| < \delta\}$.   But $T_m \in \fE - \fE$ implies that $\mu(\alpha) = \mu(-\alpha)$ whenever $\alpha \not \in \{ -1, 0, 1\}$, from which we deduce that $0 \ne \alpha \in \sigma(B)$ implies that $-\alpha \in \sigma(B)$ and $\mu(\alpha) = \mu(-\alpha)$;  in other words, $B$ is balanced.

	Since $T \in \ttt{clos}(\fE - \fE)$, by Lemma~\ref{lem3.11}, $\ttt{nul}\, (Z-I)=\ttt{nul}\, (T-I) \ge r$. Since $B$ is balanced, $\ttt{tr}(Z)=\ttt{tr}(T)=r$.
	
	\item[(b)] implies (a). \ \ \ Suppose that $T \sim B \oplus Z$, with $B$ and $Z$ as in the statement of (b).
Since $B$ is balanced, $B\in \ttt{clos}(\fE - \fE)$.  And by Lemma~\ref{lem3.10}, $Z\in \ttt{clos}(\fE - \fE)$.
Then it is clear that $T\in \ttt{clos}(\fE - \fE)$ also.
\end{enumerate}
\end{pf}

%%%%%%%%%%%%%%%%%%%%%%%%%%%%%%%%%%%%%%%%%%%

%%%%%%%%%%%%%%%%%%%%%%%%%%%%%%%%%%%%%%%%%%
%%%%%%%%%%%%%%%%%%%%%%%%%%%%%%%%%%%%%%%%%%
%%%%%%%%%%%%%%%%%%%%%%%%%%%%%%%%%%%%%%%%%%
%%%%%%%%%%%%%%%%%%%%%%%%%%%%%%%%%%%%%%%%%%
%%%%%%%%%%%%%%%%%%%%%%%%%%%%%%%%%%%%%%%%%%
%%%%%%%%%%%%%%%%%%%%%%%%%%%%%%%%%%%%%%%%%%

\section{Compact operators}

%%%%%%%%%%%%%%%%%%%%%%%%%%%%%%%%%%%%%%%%%%

\subsection{} \label{sec04.01}
Recall from Proposition~\ref{prop2.02} that
\[
\fC_\fE \subseteq \ttt{Neg}_{\ttt{invS}}(\hilb) \subseteq \ttt{Neg}_{\ttt{S}}(\hilb) \subseteq \ttt{Bal}(\hilb).\]
When $\hilb$ is finite-dimensional, the norm-closures of all of these sets coincide and  $\ttt{Bal}(\hilb)$ is closed (see Corollary~\ref{cor3.06}).
Our goal in this section is to show that the same result holds if we restrict our attention to the set of compact operators acting on an infinite-dimensional, separable Hilbert space.   That is to say, we wish to prove that
\[
\ttt{clos}\, (\fC_\fE) \cap \kofh = \ttt{Bal} (\hilb) \cap \kofh. \]

We emphasise the fact that we do not require the approximants be compact;  that is, we allow for an element $T \in \ttt{clos}\, (\fC_\fE) \cap \kofh$ (resp.  $T \in \ttt{Bal}(\hilb) \cap \kofh$) to be expressed as a limit of operators $(T_m)_m$ which are non-compact elements of $\fC_\fE$ (resp. non-compact elements of $\ttt{Bal}(\hilb)$).

\bigskip

Of course, when $K \in \kofh$ and $\lambda \in \bbC$, either $\lambda = 0$ and $K -\lambda I = K + \lambda I = K$ is not semi-Fredholm (in which case condition (iii) of Definition~\ref{defn1.07} doesn't apply), or $\lambda \ne 0$ in which case $K - \lambda I$, $K + \lambda I$ are both Fredholm of index zero, and so (iii) of Definition~\ref{defn1.07} holds automatically.
	
%%%%%%%%%%%%%%%%%%%%%%%%%%%%%%%%%%%%%%%%%%

\bigskip

The fact that every quasinilpotent operator $Q \in \bofh$ is a limit of nilpotent operators is a deep result due to Apostol and Voiculescu~\cite{ApostolVoiculescu1974}.   When $Q$ is both quasinilpotent and \emph{compact}, the fact that the approximating nilpotent operators may be chosen to be of finite rank is much simpler.   As we have been unable to locate a specific reference for this result, we have decided to include the outline of its proof.   Let $\eps > 0$ and choose a finite-rank operator $F$ such that $\norm Q - F \norm < \eps$ and $\sigma(F) \subseteq \{ z \in \bbC: |z| < \eps\}$.    That this is possible is a consequence of the fact that $Q$ is compact, combined with the upper semicontinuity of the spectrum.   Write $F \simeq F_0 \oplus 0$, where $F_0 \in \bbM_n(\bbC)$ for appropriate $n \ge 1$, and upper-triangularising $F_0$, observe that all diagonal entries have magnitude less than $\eps$.   Thus a diagonal perturbation $D_0 + F_0$ of $F_0$ of norm at most $\eps$ ($D_0$ simply represents the negative of the diagonal of $F_0$) results in a finite-rank nilpotent operator $N \simeq (D_0 + F_0) \oplus 0$ which approximates $Q$ to within $2 \eps$.

%%%%%%%%%%%%%%%%%%%%%%%%%%%%%%%%%%%%%%%%%%

\begin{prop}  \label{prop04.02}
Let $Q \in \kofh$ be quasinilpotent.   Then $Q \in \ttt{clos}\, (\fC_\fE)$.  Moreover, we can choose the approximants in $\fC_\fE$ to be nilpotent themselves.
\end{prop}

\begin{pf}
Let $\eps > 0$.   We have just seen that every compact quasinilpotent operator is a limit of finite-rank nilpotent operators, and as such, there exists a finite-rank nilpotent operator $L$ such that $\norm Q - L\norm < \eps$.   Let $\cR := \mathrm{span} \{ \mathrm{ran}\, L, \mathrm{ran}\, L^*\}$.   Then $\cR$ is finite-dimensional, and relative to the decomposition $\hilb = \cR \oplus \cR^\perp$, we may write
\[
L = \begin{bmatrix} L_0 & 0 \\ 0 & 0 \end{bmatrix}. \]
Since $L$ is nilpotent, so is $L_0$.
By Proposition~6 of~\cite{DrnovsekRadjaviRosenthal2002}, $L_0$ is a commutator of two finite-rank idempotents $E_0, F_0 \in \cB(\cR)$.    Let  $E := \begin{bmatrix} E_0 & 0 \\ 0 & 0 \end{bmatrix}$ and $F := \begin{bmatrix} F_0 & 0 \\ 0 & 0 \end{bmatrix}$.   Then $E$ and $F$ are idempotents and
\[
L = [E, F] \in \fC_\fE. \]
Since $\eps > 0$ was arbitrary, $Q \in \ttt{clos}\, (\fC_\fE)$.

\end{pf}

%%%%%%%%%%%%%%%%%%%%%%%%%%%%%%%%%%%%%%%%%%

\begin{eg} \label{eg04.03}
We temporarily digress to show that $\fC_\fE \cap \kofh$ is not closed.  Indeed, let $V\in \cB(L^2[0,1], \textup{d}x)$
be the classical Volterra operator defined by
\[
(Vf)(x)=\int_{0}^x f(t)\textup{d}t,~~f\in L^2[0,1].\]
It is well known that $V$ is compact and quasinilpotent. By Proposition~\ref{prop04.02},
we know $V\in \ttt{clos}\, (\fC_\fE)$.

A result of Kalisch (see, e.g.~\cite[Theorem~2]{Kalisch1957} or~\cite[Proposition~1]{FoiasWilliams1972}) shows that if $\alpha \in \bbC$, then $V$ and $\alpha V$ are similar if and only if $\alpha  = 1$.   In particular, $V$ is not similar to $-V$, and thus $V \not \in \fC_\fE$ by  Proposition~\ref{prop2.02}.

\end{eg}

%%%%%%%%%%%%%%%%%%%%%%%%%%%%%%%%%%%%%%%%%%

We now return to the task of extending Corollary~\ref{cor3.06} to the setting of compact operators.

\begin{prop} \label{prop04.04}
$\ttt{bal}(\hilb) \cap \kofh \subseteq  \ttt{clos}\, (\fC_\fE)$.
\end{prop}

\begin{pf}
After a moment's thought, and keeping in mind that every quasinilpotent, compact operator lies in $\ttt{clos}\, \fC_\fE$, without loss of generality, we may assume $K\in \ttt{bal}(\hilb) \cap \kofh$ is not quasinilpotent and $0$ is a cluster point of $\sigma(K)$.

\smallskip

Note that we may denote the sequence of non-zero eigenvalues of $K$ as $(\alpha_n)_n$, where $\alpha_{2k} = - \alpha_{2k-1}$, $k \ge 1$.   Since $K \in \ttt{bal}(\hilb) \cap \kofh$, it follows that $\dim\, \hilb(\{ \alpha_{2k} \}; K) = \dim\, \hilb(\{\alpha_{2k-1}\}; K)$ for all $k \ge 1$.   Define $\cM_n := \mathrm{span}\, \{ \hilb(\{ \alpha_k\}; K) \}_{k=1}^n$,  set $\hilb_\infty = \oplus_n (\cM_n \ominus \cM_{n-1})$ and $\hilb_0 = \hilb \ominus \hilb_\infty$.    Relative to the decomposition $\hilb = \hilb_\infty \oplus \hilb_0$, $K$ admits an upper triangular form
\[
K = \begin{bmatrix}
	K_{1 1} & K_{12} & K_{1 3} & \cdots & & & K_{1,0} \\
	0 & K_{2 2} & K_{2 3} & \cdots & & & K_{2, 0} \\
	0 & 0 & K_{3 3} & \cdots & & & K_{3, 0} \\
	0 & 0 & 0 & \ddots & \ddots & &  \vdots \\
	0 & 0 & 0 & 0 &  &  \cdots &    K_{0 0}
	\end{bmatrix}. \]
	
Given $n \ge 1$, let $P_n$ denote the orthogonal projection of $\hilb$ onto $\cM_n$,  let $P_0$ denote the orthogonal projection of $\hilb$ onto $\hilb_0$ and $P_\infty$ denote the orthogonal projection of $\hilb$ onto $\hilb_\infty$.

Observe that $(P_n + P_0)_n$ is an increasing sequence of projections tending strongly to the identity operator.   Since $K \in \kofh$, it follows that
\[
K = \lim_n (P_n + P_0) K (P_n + P_0). \]

Let
\[
L_{2n} := P_{2n} K P_{2n} + P_0 K P_0 =
\begin{bmatrix}
	K_{1 1} & K_{12} & K_{1 3} & \cdots & K_{1, 2n} & 0 & \cdots  & 0\\
	0 & K_{2 2} & K_{2 3} & \cdots & K_{2, 2n} & 0 & \cdots & 0 \\
	0 & 0 & \ddots & \cdots &  & 0  & \cdots  & 0 \\
	0 & 0 & 0 & \ddots & \vdots & \vdots & \cdots &   0 \\
	0 & 0 & 0 & 0 &  K_{2n, 2n} & 0 & \cdots & 0 \\
	\vdots & \vdots &  & & & & & &  \\
	0 & \cdots &  &   &    &  0 &   \cdots  & K_{0 0}
	\end{bmatrix}. \]
We may think of this as $L_{2n} = P_{2n} K P_{2n} \oplus P_0 K P_0$.   Now, $P_{2n} K P_{2n} \in \cB(\hilb_\infty)$ is a balanced, finite-rank operator.   Thus there exist finite-rank idempotents $E_{2n, \infty}$, $F_{2n, \infty} \in \cB(\hilb_\infty)$ such that $P_{2n} K P_{2n} = [E_{2n, \infty}, F_{2n, \infty}]$.

Since $\sigma(K_{00}) = \{0\}$, by Proposition~\ref{prop04.02}, there exist idempotents $E_{n, 0}$, $F_{n, 0} \in \cB(\hilb_0)$ such that $K_{0 0} = \lim_n [E_{n, 0}, F_{n, 0}]$, and each $[E_{n, 0}, F_{n, 0}]$ is nilpotent.

Then $Q_n := E_{2 n, \infty} \oplus E_{n, 0}$ and $R_n := F_{2n, \infty} \oplus F_{n, 0}$ are idempotents in $\bofh$ such that
\[
\lim_n \norm L_{2n} - [Q_n, R_n] \norm = 0. \]

Now $[Q_n, R_n] = [E_{2n, \infty}, F_{2n, \infty}] \oplus [E_{n, 0}, F_{n, 0}] \in \fC_\fE$, and therefore any operator similar to $[Q_n, R_n]$ is also in $\fC_\fE$.   Since $\sigma (P_{2n} K P_{2n}|_{\cM_{2n}}) \cap \{ 0 \} = \varnothing$, a simple argument using Rosenblum's operator (see~\cite[Corollary~3.2]{Herrero1989}) shows that
\[
[Q_n, R_n] \sim X_n := \begin{bmatrix}
	K_{1 1} & K_{12} & K_{1 3} & \cdots & K_{1, 2n} & 0 & \cdots  & K_{1, 0}\\
	0 & K_{2 2} & K_{2 3} & \cdots & K_{2, 2n} & 0 & \cdots & K_{2, 0} \\
	0 & 0 & \ddots & \cdots &  & 0  & \cdots  & \vdots \\
	0 & 0 & 0 & \ddots &  & 0 & \cdots &   \vdots  \\
	0 & 0 & 0 & 0 & K_{2n, 2n} & 0 & \cdots & K_{2n, 0} \\
	\vdots & \vdots &  & & & 0 &  \ddots  & 0 \\
	0 & \cdots &  &   &    &  0 &   \cdots  & [E_{n, 0}, F_{n, 0}]
	\end{bmatrix}. \]
Hence $X_n \in \fC_\fE$ for all $n \ge 1$.

Note that
\[
\lim_n \norm X_n - (P_{2n} + P_0) K (P_{2n} + P_0) \norm = \lim_n \norm K_{00} - [E_{n, 0}, F_{n, 0}] \norm  = 0. \]
Since $K = \lim_n (P_{2n} + P_0) K (P_{2n} + P_0)$, we conclude that
\[
K = \lim_n X_n \in \ttt{clos}\, (\fC_\fE).\]
\end{pf}

%%%%%%%%%%%%%%%%%%%%%%%%%%%%%%%%%%%%%%%%%%

We shall have reason to appeal to the next result of Herrero's more than once below.   We first recall that a \textbf{Cauchy domain} is an open set $\Omega \subseteq \bbC$ such that $\Omega$ has finitely many components, the closures of any two of which are disjoint, and the boundary $\partial(\Omega)$ consissts of a finite number of closed, rectifiable Jordan curves, any two of which are disjoint.

\begin{prop} \label{prop04.04.05} \cite[Corollary~1.6]{Herrero1989}
Let $X, Y \in \bofh$.  If $\sigma \ne \varnothing$ is a relatively closed and open subset of $\sigma(X)$, and $\Omega$ (a Cauchy domain) is a neighbourhood of $\sigma$ satisfying  $(\sigma(X)\setminus \sigma) \cap \ol{\Omega} = \varnothing$, then
\begin{itemize}
	\item{}
	$\norm X - Y \norm < \ttt{min} \{ \norm (\lambda I - X)^{-1} \norm^{-1} : \lambda \in \partial (\Omega)\}$ implies that $\sigma^\prime := \sigma(Y) \cap \Omega \ne \varnothing$; and
	\item{}
	$\dim\, \hilb(\sigma; X) = \dim\, \hilb(\sigma^\prime; Y)$.	
\end{itemize}	
\end{prop}

\smallskip

%%%%%%%%%%%%%%%%%%%%%%%%%%%%%%%%%%%%%%%%%%

\begin{prop} \label{prop04.05}
Let $K \in \kofh \cap \ttt{clos}\, (\fC_\fE)$.   Then $K \in \ttt{bal}(\hilb)$.
\end{prop}

\begin{pf}
Let $0 \ne \lambda \in \sigma(K)$.   It suffices to prove that $-\lambda \in \sigma(K)$ and that $\dim \hilb(\{ \lambda \}; K) = \dim \hilb(\{-\lambda\}; K)$.

Let $\Omega := \Omega_1 \cup \Omega_2 \cup \Omega_3$ be the \emph{disjoint} union of three open sets satisfying:
\begin{enumerate}
	\item[($\ttt{i}$)]
	$\Omega_3 = - \Omega_1$;
	\item[($\ttt{ii}$)]
	$\lambda \in \Omega_1$; (and thus $-\lambda \in \Omega_3$); and
	\item[($\ttt{iii}$)]
	$\sigma(K) \setminus \{ \lambda, -\lambda\} \subseteq \Omega_2$.
\end{enumerate}
That this is possible is clear, since $\sigma(K)$ is at most an infinite sequence of isolated points converging to zero.  In fact, we can do this while choosing $\Omega_1$ to be a disc of arbitrarily small radius centred at $\lambda$.

Let $\eps := \inf \{ \norm (\alpha I - K)^{-1} \norm^{-1} : \alpha \in \partial \Omega \} > 0$.    By Proposition~\ref{prop04.04.05}, if $T \in \bofh$ and $\norm T - K \norm < \eps$, then
\begin{enumerate}
	\item[(i)]
	$\sigma(T) \cap \Omega_1 \ne \varnothing$;
	\item[(ii)]
	$\dim( \hilb(\{ \lambda\}; K)) = \dim (\hilb(\Omega_1; T))$; and
	\item[(iii)]
	$\dim(\hilb(\Omega_3; K)) = \dim (\hilb(\Omega_3; T))$.
\end{enumerate}	
Since $K \in \ttt{clos} \, (\fC_\fE)$, there exists $X \in \fC_\fE$ such that $\norm X - K \norm < \eps$.   But $X \in \fC_\fE$ implies that $X \sim -X$, and thus $X \in \ttt{bal}(\hilb)$.    Hence
\begin{align*}
\dim(\hilb(\Omega_3; K))
	&= \dim (\hilb(\Omega_3; X)) \\
	&= \dim (\hilb(-\Omega_3; X)) \\
	&= \dim (\hilb(\Omega_1; X)) \\
	&= \dim(\hilb(\{ \lambda \}; K)) \\
	&>0.
\end{align*}
It follows that $-\lambda \in \sigma(K)$ and that
\[
\dim (\hilb(\{ -\lambda\}; K)) = \dim(\hilb(\{ \lambda \}; K)). \]
In other words, $K \in \ttt{bal}(\hilb) \cap \kofh$.
\end{pf}

%%%%%%%%%%%%%%%%%%%%%%%%%%%%%%%%%%%%%%%%%%

\begin{cor} \label{cor04.06}
\[
\kofh \cap \ttt{Bal}(\hilb) = \kofh \cap \ttt{clos}\, (\fC_\fE). \]
\end{cor}

\begin{pf}
This is an immediate consequence of Proposition~\ref{prop04.04} and Proposition~\ref{prop04.05}.
\end{pf}

%%%%%%%%%%%%%%%%%%%%%%%%%%%%%%%%%%%%%%%%%%

\begin{cor} \label{cor04.07}
\[
[\ttt{clos}\, \fC_\fE] \cap \kofh =  [\ttt{clos}\, \ttt{Neg}_{\ttt{invS}}(\hilb)] \cap \kofh =[\ttt{clos}\,  \ttt{Neg}_{\ttt{S}}(\hilb)] \cap \kofh  = \ttt{Bal}(\hilb) \cap \kofh.\]
\end{cor}

%%%%%%%%%%%%%%%%%%%%%%%%%%%%%%%%%%%%%%%%%%

\begin{eg} \label{eg4.08}
Let $K = \oplus_n \begin{bmatrix} 0 & \frac{1}{n} \\ 0 & 0 \end{bmatrix}$, so that $K$ is compact and of infinite rank.
Then  $K = [E, F]$, where $E := \oplus_n \begin{bmatrix} 1 & 0 \\ 0 & 0 \end{bmatrix}$ and $F := \oplus_n \begin{bmatrix} 1 & \frac{1}{n} \\ 0 & 0 \end{bmatrix}$ are idempotents.

Note, however, that if $E_1$ and $F_1$ are \emph{compact} idempotents, then they are necessarily of finite rank, and thus $[E_1, F_1]$ is also a finite-rank operator.   Thus $\fC_\fE \cap \kofh \ne \{ [E, F] : E, F \mbox{ are compact idempotents}\}$.
\end{eg}

%%%%%%%%%%%%%%%%%%%%%%%%%%%%%%%%%%%%%%%%%%

\bigskip

The following simple lemma will be of use to us in the next example.

\begin{lem}\label{lem4.09}
Let $\hilb_1$ and $\hilb_2$ be Hilbert spaces and $T_1=R_1\oplus 0, T_2=R_2\oplus 0\in \cB(\cH_1 \oplus \hilb_2)$, where $R_1,R_2\in \cB(\cH_1 \oplus \hilb_2)$ are invertible. If $T_1$ and $T_2$ are similar, then so are $R_1$ and $R_2$.
\end{lem}

\begin{pf}
Since $T_1$ and $T_2$ are similar,  there exists an invertible operator $S=\begin{bmatrix}A& B\\ C& D\\ \end{bmatrix} \in \cB(\hilb_1 \oplus \hilb_2)$ such that $ST_1=T_2S$.
Thus
\[ST_1=\begin{bmatrix}A& B\\ C& D\\ \end{bmatrix}\begin{bmatrix}R_1& 0 \\ 0& 0\\ \end{bmatrix}=\begin{bmatrix}AR_1& 0\\ CR_1& 0\\ \end{bmatrix} \mbox{      and      } T_2S=\begin{bmatrix}R_2& 0 \\ 0 & 0\\ \end{bmatrix}\begin{bmatrix}A& B\\ C& D\\ \end{bmatrix}=\begin{bmatrix}R_2A& R_2B\\ 0& 0\\ \end{bmatrix}.\]
From this we see that $AR_1=R_2A, CR_1=0, R_2B=0$.

As $R_1$ and $R_2$ are invertible, we conclude that $B=0$ and $C=0$. This in turn implies that $A$ is invertible in $\cB(\hilb_1)$, whence $R_1$ is similar to $R_2$ via $A$.
\end{pf}

%%%%%%%%%%%%%%%%%%%%%%%%%%%%%%%%%%%%%%%%%%

\begin{eg} \label{eg4.10}
We now produce an example of a compact operator $T$ such that $T$ is unitarily equivalent to $-T$  but it is not a commutator of idempotents in $\bofh$.  While it is based upon an example from~\cite{DrnovsekRadjaviRosenthal2002}, the extension of their result to infinite dimensions requires a surprisingly long argument.

\bigskip

Let $A_0 = \begin{bmatrix} i/2 & 1 \\ 0 & i/2 \end{bmatrix} \in \bbM_2(\bbC)$.   If $T_0 = A_0 \oplus - A_0 \in \bbM_4(\bbC)$, then by Example 10 of~\cite{DrnovsekRadjaviRosenthal2002}, $T_0$ is similar to $-T_0$ via an involution (this is trivial), but $T_0$ is not a commutator of idempotents.

Let $A = A_0 \oplus 0^{(\infty)}$, and let $T = A \oplus -A$.  Observe that $\ttt{rank}\, T = 4$.    Again, it is trivial to see that $T$ is similar to $-T$ via an involution, namely $J = \begin{bmatrix} 0 & I \\ I & 0 \end{bmatrix}$, but we claim that $T$ is not a commutator of idempotents.    (Note that $J$ is in fact a unitary involution.)

\bigskip

The proof below is an adaptation of the proof of Proposition 5 of~\cite{DrnovsekRadjaviRosenthal2002}.

\vskip 1 cm

Suppose to the contrary that $T = [E, F]$ where $E$ and $F$ are idempotents.   After conjugating by an appropriate similarity $S$, we may write
\[
S^{-1} T S = [P, Q], \]
where $P = \begin{bmatrix} I & 0 \\ 0 & 0 \end{bmatrix}$ and $Q = \begin{bmatrix} \frac{1}{2} I + B & X \\ - Y & \frac{1}{2} I - C \end{bmatrix}$ are \emph{idempotents} (but $Q$ is not necessarily a projection).     Let us assume that this decomposition of $P$ and $Q$ is relative to the decomposition $\hilb = \hilb_1 \oplus \hilb_2$ of the Hilbert space.

A calculation (which is used in~\cite{DrnovsekRadjaviRosenthal2002} and which is not hard to verify) yields that $B X = X C$ and $Y B = CY$.    It follows that $C (\ker\, X) \subseteq \ker\, X$ and that $B (\ker\, Y) \subseteq \ker\, Y$.

\bigskip

Now
\[
\ttt{rank}\, [P, Q] = 4 = \ttt{rank}\, \begin{bmatrix} 0 & X \\ Y & 0 \end{bmatrix} = \ttt{rank}\, X + \ttt{rank}\, Y. \]
Furthermore,
\[
\ttt{rank}\, (XY \oplus Y X) = [P, Q]^2 = S^{-1} T^2 S = S^{-1}  (A^2 \oplus A^2) S, \]
and therefore $\ttt{rank} \, XY + \ttt{rank} \, YX = 4$ and $\sigma(XY \oplus YX) = \sigma (T^2) = \{ -\frac{1}{4}, 0\}$.     But $\sigma(XY) \cup \{ 0\} = \sigma(YX) \cup \{ 0\}$, and thus $-\frac{1}{4} \in \sigma(XY) \cap \sigma(Y X)$, implying that $XY \ne 0 \ne Y X$.

\bigskip

Now $\ttt{rank} \, X + \ttt{rank}\, Y = 4$ from above, and neither operator is zero.   Suppose $\ttt{rank} X = 1$.   Then $\ttt{rank}\, XY \le 1$, $\ttt{rank}\, YX \le 1$ and so $\ttt{rank}\, [P, Q]^2 \le 2 \ne \rank T^2$, a contradiction.    Hence $\ttt{rank}\, X = \ttt{rank}\, Y = 2$.

Let $\cM_1 := (\ker\, X)^\perp$ and $\cM_2 := \hilb_2 \ominus \cM_1$.   Let $\cN_1 := \ker\, Y$ and $\cN_2 := \hilb_1 \ominus \cN_1$.   Relative to the decomposition $\hilb = \cN_1 \oplus \cN_2 \oplus \cM_1 \oplus \cM_2$, we may write
\[
P = \begin{bmatrix} I & 0 & 0 & 0 \\ 0 & I & 0 & 0 \\ 0 & 0 & 0 & 0 \\ 0 & 0 & 0 & 0 \end{bmatrix}, \]
and recalling that $B \cN_1 \subseteq \cN_1$ and $C \cM_2 \subseteq \cM_2$ from above,
\[
Q= \begin{bmatrix} \frac{1}{2} I + B_1 & B_2 & X_1 & 0 \\ 0 & \frac{1}{2} I + B_4 & X_3 & 0 \\ 0 & -Y_2 & \frac{1}{2} - C_1 & 0 \\ 0 & -Y_4 & -C_2 & \frac{1}{2} - C_4 \end{bmatrix}. \]

Thus, with respect to the decomposition $\hilb = \cN_1 \oplus \cM_2 \oplus \cM_1 \oplus \cN_2$, we have
\[
P = \begin{bmatrix} I & 0 & 0 & 0 \\ 0 & 0 & 0 & 0 \\ 0 & 0 & 0 & 0 \\ 0 & 0 & 0 & I \end{bmatrix}, \]
and recalling that $B \cN_1 \subseteq \cN_1$ and $C \cM_2 \subseteq \cM_2$ from above,
\[
Q= \begin{bmatrix} \frac{1}{2} I + B_1 & 0 & X_1 & B_2 \\ 0 & \frac{1}{2} - C_4 & -C_2 & -Y_4 \\ 0 & 0 & \frac{1}{2} - C_1 & -Y_2 \\ 0 & 0 & X_3  & \frac{1}{2} + B_4 \end{bmatrix}. \]

This shows that
\[
[P, Q] = \begin{bmatrix} 0 & 0 & X_1 & 0 \\ 0 & 0 & 0 &  Y_4 \\ 0 & 0 & 0 & Y_2 \\ 0 & 0 & X_3 & 0 \end{bmatrix}. \]

\bigskip

Now $\dim\, (\cM_1 \oplus \cN_2) = 4$ (since each of these spaces has dimension $2$), and
\[
[P, Q]^2 =  \begin{bmatrix} 0 & 0 & 0 & X_1 Y_2 \\ 0 & 0 & Y_4 X_3 &  0 \\ 0 & 0 & Y_2 X_3 & 0  \\ 0 & 0 & 0 & X_3 Y_2 \end{bmatrix}. \]
This is similar to $T^2$ which has four eigenvalues all equal to $-\frac{1}{4}$, and so $\sigma(Y_2 X_3) = \{ -\frac{1}{4} \} = \sigma(X_3 Y_2)$.  In particular, $X_3$ and $Y_2$ are invertible, and thus so is $\begin{bmatrix} 0 & Y_2 \\ X_3 & 0 \end{bmatrix}$.   It follows (using Rosenblum's Theorem - \cite[Corollary~3.2]{Herrero1989}) that
$[P, Q]$ is similar to
\[
R := \begin{bmatrix} 0 & 0 & 0 & 0 \\ 0 & 0 & 0 & 0 \\ 0 & 0 & 0 & Y_2 \\ 0 & 0 &  X_3 & 0 \end{bmatrix}. \]
But then by Lemma~\ref{lem4.09},
\[
\begin{bmatrix} 0 & Y_2 \\ X_3 & 0 \end{bmatrix} \mbox{ is similar to }  A_0 \oplus - A_0, \]
proving that $A_0 \oplus - A_0$ is (similar to) the commutator of the idempotents
\[
P_0 := \begin{bmatrix} 0 & 0 \\ 0 & I \end{bmatrix} \mbox{ and } Q_0 := \begin{bmatrix} \frac{1}{2} - C_1 & - Y_2 \\ X_3 & \frac{1}{2} I + B_4 \end{bmatrix}. \]
\end{eg}

%%%%%%%%%%%%%%%%%%%%%%%%%%%%%%%%%%%%%%%%%%

\begin{notation}
Given  $\varnothing \ne L \subseteq \bbC$ and $\eps > 0$, we define $L_\eps := \{ z \in \bbC: \ttt{dist}\, (z, L) < \eps \}$.

\bigskip

Recall that if $\Delta := \{ A \subseteq \bbC: A \mbox{ is compact}\}$, then the \textbf{Hausdorff metric} on $\Delta$ is the metric defined by
\[
 d_H(A, B)  :=  \ttt{max}\  ( \ttt{max}_{a \in A} \ttt{dist}(a, B),  \ttt{max}_{b \in B} \ttt{dist} (b, A)). \]
\end{notation}		
%Since the way normality is  used is via the fact that for normal operators
% the spectral radius is the same as the norm,  we could extend $mutatis~mutandis$ this result from $\cB(\cH)$ to a unital $C^*$-algebra $\cA$.

%%%%%%%%%%%%%%%%%%%%%%%%%%%%%%%%%%%%%%%%%%

\begin{thm} \label{thm4.11}
Let $K \in \ttt{clos}\, (\fE - \fE) \cap \kofh$, and write
\[
K = \begin{bmatrix} K_1 & K_2 \\ 0 & K_4 \end{bmatrix} \]
relative to the decomposition $\hilb = \hilb(\{ -1, 1\}; K) \oplus (\hilb(\{ -1, 1\}; K))^\perp$.     By considering $-K$ instead of $K$ if necessary, we may assume without loss of generality that
\[
\ttt{tr}\, (K_1) \ge 0. \]
Then
\begin{enumerate}
	\item[(i)]
	$\ttt{nul} \, (K - I) \ge \ttt{tr}\, K_1$; and
	\item[(ii)]
	$K_4$ is balanced.
\end{enumerate}
\end{thm}

\noindent{\textbf{Remarks.}}\ \ \
An equivalent formulation of (ii) is that if $0 \ne \alpha \in \sigma(K) \setminus \{ -1, 1\}$, then $-\alpha \in \sigma(K)$ and $\mu(\alpha) = \mu(-\alpha)$.
Also, if $\sigma(K) \cap \{ -1, 1\} = \varnothing$, then $K_1$, $K_2$ above are absent and $K = K_4$ in the argument below, meaning that one is only required to prove that $K_4$ is balanced.  As such, the first half of the proof (regarding $K_1$) only applies if $\sigma(K) \cap \{ -1, 1 \} \ne \varnothing$.

\smallskip

\begin{pf}
Since $K$ is compact, we know that $0 \ne \alpha \in \sigma(K)$ implies that $\alpha$ is isolated, and thus there exists $\delta > 0$ such that if $G_{-1} := \{ -1\}_\delta$, $G_1 := \{ 1 \}_\delta$ and $G_0 := (\sigma(K)\setminus \{ -1, 1\})_\delta$, then (using $\sqcup$ to denote the \emph{disjoint} union of sets)
\[
\sigma(K) \subseteq G_{-1} \sqcup G_{1} \sqcup G_0. \]
Let $(T_m)_m$ be a sequence in $\fE - \fE$ such that $K  = \lim_m T_m$.   By Proposition~\ref{prop04.04.05}, there exists $m_0 \ge 1$ such that $m \ge m_0$ implies that
\begin{enumerate}
	\item[(a)]
	$\sigma(T_m) \subseteq G_{-1} \sqcup G_1  \sqcup G_0$;
	\item[(b)]
	$\sigma(T_m) \cap G_{-1} \ne \varnothing$, $\sigma(T_m) \cap G_{1} \ne \varnothing$,  $\sigma(T_m) \cap G_0 \ne \varnothing$; and
	\item[(c)]
	$\dim\, \hilb(G_1; T) = \dim\, \hilb(G_1; K) < \infty$ and $\dim\, \hilb(G_{-1}; T) = \dim\, \hilb(G_{-1}; K) < \infty$.
\end{enumerate}
Clearly $\ttt{tr}\, K_1 = \dim\, \hilb(G_1; K) - \dim\, \hilb(G_{-1}; K)$.

Meanwhile,  if we write
\[
T_m = \begin{bmatrix} T_1^{(m)} & T_2^{(m)} \\ 0 & T_4^{(m)} \end{bmatrix} \]
relative to the decomposition $\hilb = \hilb(G_1 \sqcup G_{-1}; T_m) \oplus (\hilb(G_1 \sqcup G_{-1}; T_m))^\perp$, then $T_1^{(m)}$ acts on a finite-dimensional space and $\ttt{tr}(T_1^{(m)}) = \sum \{ \mu(\beta) \beta:  \beta \in \sigma(T_m) \cap (G_1 \sqcup G_{-1}) \}$.

Since $\sigma(-T_1^{(m)}) \cap \sigma(T_4^{(m)}) = \varnothing$, by Lemma~2 of Hartwig and Putcha~\cite{HartwigPutcha1990} (the reader should be aware that there is a typographical error in the statement of their Lemma -- the correct hypothesis there should be that $\sigma(-P) \cap \sigma(Q) = \varnothing$, as they require in their proof), $T_1^{(m)} \in \fE - \fE$.  From this and their characterisation of $\fE - \fE$ (i.e. the eigenvalues of $T_1^{(m)}$ which are different from $0, -1$ and $1$ come in pairs when counted with algebraic multiplicity) and Lemma~\ref{lem3.11}, it follows that
\begin{align*}
\ttt{tr}\, K_1
	&= \dim\, \hilb(G_1; K) - \dim\, \hilb(G_{-1}; K) \\
	&= \dim\, \hilb(G_1; T_m) - \dim \, \hilb(G_{-1}; T_m) \\
	&= \dim\, \hilb(G_1; T_1^{(m)}) - \dim\, \hilb(G_{-1}; T_1^{(m)}) \\
	&= \dim\, \hilb(\{ 1 \}; T_1^{(m)}) - \dim\, \hilb(\{-1\}; T_1^{(m)}) \\
	&= \ttt{tr}(T_1^{(m)}) \\
	&\le \ttt{nul} (T_1^{(m)} - I) \\
	&= \ttt{nul}(T_m - I).
\end{align*}	
	
But $K = \lim_m T_m$, and thus $\ttt{nul} \, (K-I) \ge \ttt{nul}\, (T_m - I)$, whence $\ttt{nul}\, (K-I) \ge \ttt{tr}\, K_1$, proving (i) above.

\bigskip

There remains to show that $K_4$ is balanced.  Obviously it suffices to consider the case where $K \ne 0$.   Fix a strictly decreasing sequence $(\delta_n)_n$ of strictly positive real numbers satisfying:
\begin{itemize}
	\item[(i)]
	$\delta_2 < 1 < \delta_1 < \delta_0 := \norm K \norm + 1$;
	\item[(ii)]
	$\alpha \in \sigma(K)$ and $\delta_2 < |\alpha| < \delta_1$ implies that $|\alpha| = 1$;
	\item[(iii)]
	$\alpha \in \sigma(K) $ implies that $|\alpha| \not \in \{ \delta_n\}_n$; and
	\item[(iv)]
	$\lim_n \delta_n = 0$.
\end{itemize}
Since $\sigma(K)$ is a sequence converging to $0$, this is easy to do.

\bigskip

Given $0 < r < s$, let $A_{r, s} := \{ z \in \bbC: r < |z| < s\}$ be the open annulus centred at $0$ of inner radius $r$ and outer radius $s$.    Abbreviate this to $\Omega_n := A(\delta_{n+1}, \delta_n)$, $n \ge 0$.

For any $n \ge 0$, it is clear that $\dim \, \hilb(\Omega_n; K) < \infty$.  Also, if $\{ -1, 1 \} \cap \sigma(K) \ne \varnothing$, then $\{ - 1, 1 \} \cap \sigma(K) \subseteq \Omega_1$.

Let $\eps > 0$ and choose $N \ge 1$ such that $\delta_N < \eps$.    Let $\Omega := \cup_{n=0}^N \Omega_n$, and let $\Gamma := \{ z \in \bbC: |z| < \delta_{N+1} \}$.   Then $\Omega, \Gamma$ are open and disjoint, and $\sigma(K) \subseteq \Omega \cup \Gamma$.  In fact, for $0 \le j \le N$, $\sigma(K) \cap \Omega_n$ is a finite set (including multiplicity).   Define
\[
\eta_1 := \frac{1}{2} \min \{ |\alpha - \beta| : \alpha, \beta \in \Omega \cap \sigma(K), \alpha \ne \beta \}, \]
\[
\eta_2 := \min \{ \ttt{dist} (\alpha, \partial \Omega_n) : \alpha \in \Omega_n \cap \sigma(K), 0 \le n \le N\}, \]
and choose $0 < \eta < \min (\eta_1, \eta_2)$.
\bigskip

Then we may find open sets $\Omega_j^{(n)}$, $1 \le j \le r_n$, $0 \le n \le N$ such that
\begin{enumerate}
	\item[(\textsc{i})]
	$\sigma(K) \cap \Omega_j^{(n)}$ contains exactly one element for all $1 \le j \le r_n$; $0 \le n \le N$ (though possibly with algebraic multiplicity greater than one);
	\item[(\textsc{ii})]
	$\sigma(K) \cap \Omega = \cup_{n=0}^N \cup_{j=1}^{r_n} \sigma(K) \cap \Omega_j^{(n)}$;
	\item[(\textsc{iii})]
	$\ttt{diam} \, \Omega_{j}^{(n)} < \eta$ for all $1 \le j \le r_n$, $0 \le n \le N$; and
	\item[(\textsc{iv})]
	$\ttt{diam}\, \Omega_j^{(n)} = \ttt{diam}\, \Omega_k^{(n)}$ for all $1 \le j, k \le r_n$, $0 \le n \le N$.
\end{enumerate}
(In essence, we take a ball of radius $\eta$ around each $\alpha \in \sigma(K) \cap \Omega$, and observe that  $\eta < \eta_1$ ensures that each such ball only contains one element of $\sigma(K)$, and that no two such balls intersect.  Furthermore, if $\alpha \in \Omega_n$ for some $0 \le n \le N$, then $\eta < \eta_2$ implies that the entire ball of radius $\eta$ centred at $\alpha$ is contained in that $\Omega_n$.)

\bigskip

By the upper-semicontinuity of the spectrum and (an induction argument using~Proposition~\ref{prop04.04.05}), there exists $\zeta > 0$ such that if $T \in \bofh$ and $\norm T - K \norm < \zeta$, then
\begin{enumerate}
	\item[(a)]
	$\sigma(T) \subseteq \left( \cup_{0 \le n \le N}\  \left(\cup_{1 \le j \le r_n} \Omega_j^{(n)} \right) \right) \cup \Gamma$; and
	\item[(b)]
	$\dim\, \hilb(\Omega_j^{(n)}; T) = \dim\, \hilb(\Omega_j^{(n)};K) < \infty$, $1 \le j \le r_n$, $0 \le n \le N$.
\end{enumerate}
In particular, since $K \in \ttt{clos}\, (\fE -\fE)$, we may assume that $T \in \fE -\fE$ and that $\norm T - K \norm < \zeta$.

\bigskip

The fact that $\sigma(T) \subseteq \cup_{n=0}^N \Omega_n \cup \Gamma$ and that these sets are open and disjoint ensures that relative to the decomposition $\hilb = \oplus_{n=0}^N  \hilb(\Omega_n; T) \oplus \hilb(\Gamma;T)$, we may write $T$ as an upper-triangular operator matrix
\[
T = [T_{i,j}]. \]
Note that for all $0 \le n \le N$,   $\sigma(T_{n, n}) \subseteq \Omega_n$, and thus if $0 \le i \ne j \le N$, then  $\sigma(-T_{i,i}) \cap \sigma(T_{j,j}) = \varnothing$.   Furthermore, $\Omega_n \cap \Gamma = \varnothing$ for all $0 \le n \le N$, and $\sigma(T_{N+1 ,N+1}) \subseteq \Gamma$.   From this we conclude that
\[
T \sim \ttt{diag} (T_{0,0}, T_{1,1}, T_{2,2}, \ldots, T_{N, N}, T_{N+1, N+1}). \]
By the Hartwig-Putcha Theorem~\cite{HartwigPutcha1990}, since $T \in \fE-\fE$, it follows that each $T_{n, n} \in \fE-\fE$, $0 \le n \le N+1$.    But $T_{n, n} \in \cB(\hilb(\Omega_n;T))$, and $\dim\, \hilb(\Omega_n; T) < \infty$, $0 \le n \le N$.   By the Hartwig-Putcha characterisation of $\ttt{doi}$s in the finite-dimensional setting, each $\sigma(T_{n, n})$ is balanced, $0 \le n \le N, \ n \ne 1$; and $\sigma(T_{1,1}) \setminus \{-1, 1\}$ is balanced.

\bigskip

Let $\alpha \in \sigma(K) \cap \Omega$, $\alpha \not \in \{ -1, 1\}$, and choose $0 \le n \le N$, $1 \le j \le r_n$ such that $\alpha \in \Omega_j^{(n)}$.    From (b) above,
\[
\dim\, \hilb(\Omega_j^{(n)}; T) = \dim\,\hilb(\Omega_j^{(n)};K) = \mu(\alpha), \]
the algebraic multiplicity of $\alpha$ in $\sigma(K)$.  Since $T$ is balanced,
\[
\dim\, \hilb(\Omega_j^{(n)}; T) = \dim\, \hilb(-\Omega_j^{(n)}; T). \]

Now from condition (a) above, the Hausdorff metric
\[
\ttt{dist}_H (\sigma(T) \cap \Omega, \sigma(K) \cap \Omega) < \eta, \]
and thus there exists $\beta \in \sigma(K)$ such that $|\beta -  (-\alpha)| < 2 \eta$.     But $\eta > 0$ can be chosen arbitrarily small, implying that $-\alpha \in \sigma(K) \cap \Omega_n$.  Since $\alpha \ne 0$, $-\alpha \in \Omega_i^{(n)}$ for some $1 \le i \ne j \le r_n$ and hence $\Omega_i^{(n)} = - \Omega_j^{(n)}$.

But then
\begin{align*}
\dim\, \hilb(\Omega_j^{(n)};K)
	&= \dim\, \hilb(\Omega_j^{(n)}; T) \\
	&= \dim\, \hilb(-\Omega_j^{(n)}; T) \\
	&= \dim\, \hilb(\Omega_i^{(n)}; T) \\
	&= \dim\, \hilb(\Omega_i^{(n)}; K) \\
	&= \dim\, \hilb(-\Omega_j^{(n)}; K),
\end{align*}	
which implies that $(\sigma(K) \setminus \{ -1, 1\}) \cap \Omega$ is balanced.

\bigskip

Recall that $\Omega = \cup_{0\le n \le N} \Omega_n$, and that the only condition on $N\ge 1$ was that we must have $\delta_N < \eps$.   In particular, we can choose $N$ arbitrarily large, and from this we conclude that if $\alpha \in \sigma(K) \setminus \{ -1, 1\}$, then $-\alpha \in \sigma(K)$ and $\mu(\alpha) = \mu(-\alpha)$.   This completes the proof.

\end{pf}
	
%%%%%%%%%%%%%%%%%%%%%%%%%%%%%%%%%%%%%%%%%%

\begin{thm} \label{thm4.12}
Let $K \in \kofh$ and write
\[
K = \begin{bmatrix} K_1 & K_2 \\ 0 & K_4 \end{bmatrix} \]
relative to the decomposition $\hilb = \hilb(\{ -1, 1 \}; K) \oplus (\hilb(\{-1, 1\}; K))^\perp$.     Without loss of generality (by considering $-K$ instead of $K$ if necessary), we may suppose that $\ttt{tr}\, (K_1) \ge 0$.   The following are equivalent:
\begin{enumerate}
	\item[(a)]
	$K \in \ttt{clos}(\fE - \fE)$.
	\item[(b)]
	$\ttt{nul}\, (K_1 - I) \ge \ttt{tr}\, (K_1)$ and  $\sigma(K_4)$ is balanced.
	\end{enumerate}	
\end{thm}

\begin{pf}
\begin{enumerate}
	\item[(a)] implies (b). \ \ \
	This is Theorem~\ref{thm4.11}.
	\item[(b)] implies (a). \ \ \ If $\sigma(K) \cap \{ -1, 1\} = \varnothing$, then $K = K_4$ and there remains only to show that $K_4$ is balanced.   Otherwise, observe that since $\sigma(K_1) \subseteq \{ -1, 1\}$ and $\sigma(K_4) \subseteq \sigma(K) \setminus \{ -1, 1\}$, we have that $\sigma(-K_1) \cap \sigma(K_4) = \varnothing$, whence $K \sim \begin{bmatrix} K_1 & 0 \\ 0 & K_4 \end{bmatrix}$.
Since $\dim\, \hilb(\{ -1, 1\}; K) < \infty$, $K_1$ acts on a finite-dimensional space, and so we may apply Lemma~\ref{lem3.10} to conclude that $K_1 \in \ttt{clos}\, (\fE-\fE)$.

Since $K_4$ is balanced, $K_4 \in \ttt{clos}\, (\fE - \fE)$ by Proposition~\ref{prop04.04}.

It is now clear that $K_1 \oplus K_4 \in \ttt{clos}\, (\fE-\fE)$.   Since $K$ is similar to $K_1 \oplus K_4$ and since $\fE - \fE$ is invariant under conjugation by invertible elements, $K \in \ttt{clos}\, (\fE - \fE)$.
\end{enumerate}
\end{pf}

%%%%%%%%%%%%%%%%%%%%%%%%%%%%%%%%%%%%%%%%%%
%%%%%%%%%%%%%%%%%%%%%%%%%%%%%%%%%%%%%%%%%%
%%%%%%%%%%%%%%%%%%%%%%%%%%%%%%%%%%%%%%%%%%
%%%%%%%%%%%%%%%%%%%%%%%%%%%%%%%%%%%%%%%%%%
%%%%%%%%%%%%%%%%%%%%%%%%%%%%%%%%%%%%%%%%%%
%%%%%%%%%%%%%%%%%%%%%%%%%%%%%%%%%%%%%%%%%%

\section{Commutators and differences of orthogonal projections} \label{sec5}

\subsection{} \label{sec5.01}
Our goal in this section is to describe the sets $\ttt{clos}\, (\fC_\fP)$ and $\ttt{clos}\, (\fP - \fP)$.   We are aided by the fact that the sets $\fC_\fP$ and $\fP - \fP$ have been completely characterised by Li~\cite{Li2004} and Davis~\cite{Davis1958} respectively.  We begin with $\ttt{clos}\, (\fC_\fP)$.

%%%%%%%%%%%%%%%%%%%%%%%%%%%%%%%%%%%%%%%%%%

\begin{thm} \label{thm5.02} \emph{\textbf{[Li.]}}\ \ \
Let $\hilb$ be a complex Hilbert space.   An operator $T \in \bofh$ is a commutator of two orthogonal projections if and only if
\begin{enumerate}
	\item[(a)]
	$T^* = - T$;
	\item[(b)]
	$\norm T \norm \le \frac{1}{2}$; and
	\item[(c)]
	$T \simeq T^*$.
\end{enumerate}
\end{thm}

%%%%%%%%%%%%%%%%%%%%%%%%%%%%%%%%%%%%%%%%%%

\subsection{} \label{sec5.03}
It is worth observing that when $n := \dim \hilb < \infty$, both $\fC_\fP$ and $\fP - \fP$ are norm-closed.   Indeed, if $(P_n)_n$, $(Q_n)_n$ are two sequences of orthogonal projections in $\cB(\bbC^n)$, then the fact that the closed unit ball of $\cB(\bbC^n)$ is compact can be used to prove that there exists a strictly increasing sequence $(n_k)_k$ of positive integers such that $P:= \lim_k P_{n_k}$ and $Q := \lim_k Q_{n_k}$ both exist.   Clearly both $P$ and $Q$ are orthogonal projections, and thus if $T = \lim_n [P_n, Q_n]$, we conclude that $T = [P, Q]$, while  if $R = \lim_n (P_n - Q_n)$, then $R = P- Q$.

%%%%%%%%%%%%%%%%%%%%%%%%%%%%%%%%%%%%%%%%%%

\subsection{}  \label{sec5.04}
We remark that unlike the situation with idempotents, $\fC_\fP \not \subseteq \fP - \fP$.   For example, if $P = \begin{bmatrix} 1 & 0 \\ 0 & 0 \end{bmatrix}$ and $Q = \begin{bmatrix} \frac{1}{2} & \frac{1}{2} \\ \frac{1}{2} & \frac{1}{2} \end{bmatrix}$, then
\[
[P, Q] = \begin{bmatrix} 0 & \frac{1}{2} \\ -\frac{1}{2} & 0 \end{bmatrix} \ne [P, Q]^*, \]
while any difference of orthogonal projections is clearly self-adjoint.  Indeed, condition (a) of Li's Theorem above implies that $\fC_\fP \cap (\fP - \fP) = \{ 0\}$.

%%%%%%%%%%%%%%%%%%%%%%%%%%%%%%%%%%%%%%%%%%

\begin{lem} \label{lem5.05}
Let $\hilb$ be a complex Hilbert space and $K = K^* \in \bofh$.   Suppose furthermore that $K \simeq_a -K$.    Given $\eps > 0$, there exists an operator  $L_\eps = L_\eps^* \in \bofh$ satisfying
\begin{enumerate}
	\item[(a)]
	$\norm L_\eps \norm \le \norm K \norm$;
	\item[(b)]
	$\norm L_\eps - K \norm <  2 \eps$; and
	\item[(c)]
	$L_\eps \simeq - L_\eps$.
\end{enumerate}
\end{lem}
	
\begin{pf}
Clearly it suffices to consider the case where $\norm K \norm = 1$.

\smallskip

By hypothesis, $K$ is self-adjoint (hence normal) and $K$ is approximately unitarily equivalent to $-K$.   It follows from the Weyl-von Neumann-Berg Theorem (see, e.g.~\cite[Theorem~II.4.4]{Davidson1996}) that $\sigma(K) = - \sigma(K)$, and if $\alpha \in \sigma(K)$ is isolated, then
\[
\ttt{nul}\, (K - \alpha I) = \ttt{nul}\, (K+ \alpha I). \]

Thus, if $0< \eps < 1$, then $\sigma(K) \cap [\eps, 1] = - (\sigma(K) \cap [-1, -\eps])$, including the multiplicity (finite or infinite) of isolated eigenvalues.  Let $\cM_\eps^+ := \hilb(\sigma(K) \cap [\eps, 1]; K)$ and $\cM_\eps^- := \hilb(\sigma(K) \cap [-1, -\eps]; K)$.   Set $\cN_\eps := \hilb \ominus (\cM_\eps^+ \oplus \cM_\eps^-)$.

Again, by the Weyl-von Neumann-Berg Theorem,
\[
K|_{\cM_\eps^+} \simeq_a - K|_{\cM_\eps^-}. \]
%Denote by $P_\eps^+$ (resp. $P_\eps^-$) the orthogonal projection of $\hilb$ onto $\cM_\eps^+$ (resp. onto $\cM_\eps^-$), and let $Q_\eps$ denote the orthogonal projection of $\hilb$ onto $\cN_\eps$.

\smallskip

Relative to $\hilb = \cM_\eps^- \oplus \cN_\eps \oplus \cM_\eps^+$, we may write
\[
K= \begin{bmatrix} K_\eps^- & & \\  & K_\eps^\circ & \\  & & K_\eps^+ \end{bmatrix}. \]
From above, $K_\eps^- \simeq_a - K_\eps^+$, and so we can find a unitary operator $V$ such that
\[
\norm V^* (-K_\eps^+) V - K_\eps^-\norm < \eps. \]
Also, $\norm K_\eps^\circ \norm \le \eps$.   Define (with respect to the same decomposition of $\hilb$) the operator
\[
L_\eps :=  \begin{bmatrix} V^*(- K_\eps^+) V & & \\  & 0 & \\  & & K_\eps^+ \end{bmatrix}. \]

Clearly $L_\eps \simeq -L_\eps$, and $\norm L_\eps - K \norm \le \eps < 2 \eps$.    Note also that $\norm L_\eps \norm = \norm K_\eps^+\norm \le \norm K \norm$.
\end{pf}

%%%%%%%%%%%%%%%%%%%%%%%%%%%%%%%%%%%%%%%%%%

\begin{rem} \label{rem5.06}
We note that from the construction of $L_\eps$ above, if $\norm K \norm$ is not an eigenvalue of $K$, then neither $ - \norm K \norm$  nor $\norm K \norm$ are eigenvalues of $L_\eps$.
\end{rem}

%%%%%%%%%%%%%%%%%%%%%%%%%%%%%%%%%%%%%%%%%%

\smallskip

We are now in a position to characterise $\ttt{clos}\, \fC_\fP$.

\smallskip

\begin{thm}  \label{thm5.07}
Let $\hilb$ be a complex Hilbert space.   An operator $T \in \bofh$ is a limit of commutators of projections if and only if it satisfies the following three conditions:
\begin{enumerate}
	\item[(a)]
	$T^* = - T$;
	\item[(b)]
	$\norm T \norm \le \frac{1}{2}$; and
	\item[(c)]
	$T$ is approximately unitarily equivalent to $T^*$.
\end{enumerate}
\end{thm}
	
\begin{pf}
Suppose that $T$ satisfies the above three conditions.  Let $K := i T$.   Then $K^* = (i T)^* = -i (-T) = i T = K$, $\norm K \norm = \norm T \norm \le \frac{1}{2}$, and if $T^* = \lim_n U_n^* T U_n$, where each $U_n$ is unitary, $n \ge 1$, then
\[
K = K^* = -i T^* = -i  \lim_n U_n^* T U_n = \lim_n U_n^* (- i T) U_n = \lim_n U_n (-K) U_n. \]
That is, $K$ is self-adjoint and approximately unitarily equivalent to $-K$.

Let $\eps > 0$, and using Lemma~\ref{lem5.05}, we may choose $L_\eps = L_\eps^*$ such that
\begin{enumerate}
	\item[(a)]
	$\norm L_\eps \norm \le \norm K \norm \le \frac{1}{2}$;
	\item[(b)]
	$\norm L_\eps - K \norm < 2 \eps$; and
	\item[(c)]
	$L_\eps \simeq -L_\eps$.
\end{enumerate}	

By Theorem~\ref{thm5.02},  $T_\eps := -i L_\eps \in \fC_\fP$ and clearly $\norm T - T_\eps \norm < 2 \eps$.  It follows that $T \in \ttt{clos}\, (\fC_\fP)$.

The reverse containment is straightforward and is left to the reader.
\end{pf}

%%%%%%%%%%%%%%%%%%%%%%%%%%%%%%%%%%%%%%%%%%

Our next goal is to classify the closure $\ttt{clos}\, (\fP - \fP)$ of the set $\fP - \fP$ of differences of projections in $\bofh$, where $\hilb$ is an arbitrary complex Hilbert space.   We first recall the Theorem of Davis~\cite[Theorem~6.1]{Davis1958}.

\begin{thm} \label{thm5.08} \emph{\textbf{[Davis.]}} \ \ \
Let $H \in \bofh$ be a self-adjoint operator of norm at most one, and define $\hilb_0 := (\ker (H + I) \oplus \ker H \oplus \ker (H-I))^\perp$.    The following are equivalent:
\begin{enumerate}
	\item[(a)]
	$H \in \fP - \fP$.
	\item[(b)]
	$H_0 \simeq -H_0$, where $H_0 := H|_{\hilb_0}$.
\end{enumerate}	
\end{thm}

%%%%%%%%%%%%%%%%%%%%%%%%%%%%%%%%%%%%%%%%%%

We require  a couple of standard results;  the first is due to Newburgh~\cite{Newburgh1951} (alternatively, see~\cite[Problem~105]{Halmos1982}).

\begin{thm} \emph{\textbf{[Newburgh.]}}\ \ \  \label{thm5.09}
Let $(M_k)_k$ be a sequence of normal operators on $\hilb$ which converge in norm to $M \in \bofh$.    Then $(\sigma(M_k))_k$ converges to $\sigma(M)$ in the Hausdorff metric.
\end{thm}

%%%%%%%%%%%%%%%%%%%%%%%%%%%%%%%%%%%%%%%%%%

\bigskip

The \textbf{essential spectrum} of an operator $T \in \bofh$ is the spectrum $\sigma_e(T) := \sigma(\pi(T))$, where $\pi: \bofh \to \bofh/\kofh$ is the canonical quotient map.   When $T$ is a normal operator, the relationship between $\sigma(T)$ and $\sigma_e(T)$ is particularly simple, and is given by the next result~\cite[Proposition~4.6]{Conway1990}.

\begin{prop}\label{prop5.10}
If $N \in \bofh$ is a normal operator, then
\[
\sigma(N) \setminus \sigma_e(N)=\{\lambda\in \sigma(N): \lambda~\textup{is}~\textup{an~isolated~eigenvalue~of~finite~multiplicity~of}~N\}.\]
\end{prop}

%%%%%%%%%%%%%%%%%%%%%%%%%%%%%%%%%%%%%%%%%%

\smallskip

We now have all the tools we need to characterise  the set $\ttt{clos} \, (\fP - \fP)$.

\smallskip

\begin{thm}  \label{thm5.11}
Let $\hilb$ be a complex Hilbert space.   An operator $H \in \bofh$ is a limit of differences of projections if and only if it satisfies the following two conditions:
\begin{enumerate}
	\item[(a)]
	$-I \le H \le I$; and
	\item[(b)]
	if $\cN := (\ker (H^2 - I))^\perp$, and $H_1 := H|_{\cN}$, then $H_1 \simeq_a - H_1$.
\end{enumerate}
\end{thm}
	
\begin{pf}
	Suppose that $K_n \in \fP - \fP$, $n \ge 1$ and that $H = \lim_n K_n$.    Since $K_n = K_n^*$ for all $n\ge 1$, we have that $H = H^*$.   Also, 	
	since $- I \le K_n \le I$ for all $n\ge 1$, we have that $-I \le H \le I$.

	It is now a consequence of the Weyl-von Neumann-Berg Theorem~\cite[Theorem~II.4.4]{Davidson1996} and Proposition~\ref{prop5.10} that to prove that $H_1 \simeq_a -H_1$, it suffices to show that $\sigma(H_1) =-\sigma(H_1)$, and if $\alpha \in \sigma(H_1)$ satisfying $| \alpha | < 1$ is an isolated eigenvalue of finite multiplicity $\mu(\alpha)$, then $-\alpha \in \sigma(H_1)$ is an isolated eigenvalue of the same multiplicity.   (Note that $\sigma(H_1) = -\sigma(H_1)$ implies that $\sigma(H_1) \subseteq [-1, 1]$ is symmetric about the origin.   Furthermore, by definition of $\cN$, $1 \in \sigma(H_1)$ if and only if $1$ is a limit of a sequence $(\alpha_n)_n$ in $\sigma(H_1) \cap (-1, 1)$, in which case $-\alpha_n \in \sigma(H_1)$ for all $n \ge 1$, and thus $-1 = \lim_n -\alpha_n \in \sigma(H_1)$ as well.)
	
	Recall that $H = \lim_n K_n$.    By Newburgh's Theorem~\ref{thm5.09},
	\[
	\lim_n d_H(\sigma(K_n), \sigma(H)) = 0. \]
%	\[
%	\lim_n d_H(\sigma_e(K_n), \sigma_e(H)) = 0. \]
	Hence, if $\alpha \in \sigma(H)$ with $| \alpha | < 1$, there exists $\alpha_n \in \sigma(K_n)$ with $\lim_n \alpha_n = \alpha$.   Of course, this in turn implies that there exists $N_1 \ge 1$ such that $n \ge N_1$ implies that $| \alpha_n | < 1$.
	
	But $K_n \in \fP - \fP$, and thus  by Davis' Theorem~\ref{thm5.08}, $-\alpha_n \in \sigma(K_n)$.    By Newburgh's Theorem~\ref{thm5.09} (once again), $-\alpha  = \lim_n -\alpha_n \in \sigma(H)$.  Thus $\sigma(H_1) \cap (-1, 1) = - \sigma(H_1) \cap (-1, 1)$, and from this we see that $\sigma(H_1) = - \sigma(H_1)$.
	
	Suppose next that $\alpha \in \sigma(H_1) \cap (-1, 1)$ is an isolated eigenvalue of finite multiplicity $\mu(\alpha)$.   Since $\sigma(H_1) = - \sigma(H_1)$, $-\alpha \in \sigma(H_1)$ is an isolated eigenvalue of multiplicity $\mu(-\alpha)$.   Fix $0 < \eps < \min( |\alpha|,  1 - |\alpha|)$ such that
	\[
	(\alpha - \eps, \alpha + \eps) \cap \sigma(H_1) = \{ \alpha \}. \]
	By the symmetry of $\sigma(H_1) \subseteq \bbR$ about the origin, we have that
	\[
	(-\alpha - \eps, -\alpha + \eps) \cap \sigma(H_1) = \{ -\alpha\}. \]
	Furthermore,
	\[
	\mu(\alpha) = \dim \, \hilb( (\alpha-\eps, \alpha+\eps) \cap \sigma(H_1); H_1). \]
	By Proposition~\ref{prop04.04.05}, there exists $N_2 \in \bbN$ such that $n \ge N_2$ implies that
	\[
	\dim\ \hilb((\alpha - \eps, \alpha+ \eps) \cap \sigma (K_n); K_n)  = \mu(\alpha), \mbox{ and }\]
	\[
	\dim \, \hilb( (-\alpha-\eps, -\alpha+\eps) \cap \sigma(K_n); K_n) = \mu(-\alpha). \]
	
	\smallskip
	
	The fact that each $K_n \in \fP - \fP$ implies (by Davis' Theorem)  that
	\[
	\dim \, \hilb( (\alpha-\eps, \alpha+\eps) \cap \sigma(K_n); K_n) = \dim \, \hilb( (-\alpha-\eps, -\alpha+\eps) \cap \sigma(K_n); K_n), \]
	whence
	\[
	\mu(\alpha) = \mu(-\alpha), \]
	completing the proof of the fact that $H \in \fP - \fP$ implies both conditions (a) and (b).
	\bigskip

	Conversely, suppose that $H$ satisfies (a) and (b) above.   Of course, $-I < H_1 < I$ is an hermitian operator.  Relative
	to $\hilb = \cN^\perp \oplus \cN$, we may write
	\[
	H = \begin{bmatrix} H^\circ & \\ & H_1 \end{bmatrix}. \]
	(Note that $H^\circ$ is hermitian with $\sigma(H^\circ) \subseteq \{ -1, 1\}$.)   Let $\eps > 0$.   Using
	Lemma~\ref{lem5.05}, we can find an hermitian operator  $L_\eps$ such that
	\begin{enumerate}
		\item[(a)]
		$\norm L_\eps \norm \le \norm H_1 \norm \le 1$;
		\item[(b)]
		$\norm L_\eps - H_1 \norm < 2 \eps$; and
		\item[(c)]
		$L_\eps \simeq -L_\eps$.
	\end{enumerate}	

	Let $H_\eps := \begin{bmatrix} H^\circ & \\ & L_\eps \end{bmatrix}$.  As noted in Remark~\ref{rem5.06} (and
	keeping in mind that $\norm L_\eps \norm \le \norm H_1 \norm$), since $1$ is not an eigenvalue of $H_1$,
	neither $-1$ nor $1$ are eigenvalues of $L_\eps$.      It is clear that  $\norm H_\eps - H \norm < 2 \eps$
	and $H_\eps \in \fP - \fP$ by Theorem~\ref{thm5.08}.

	Thus $H \in \ttt{clos}\, \fP - \fP$.
\end{pf}

\subsection*{Acknowledgements}
Laurent W. Marcoux's research was supported in part by\linebreak NSERC (Canada). Yuanhang Zhang's research was supported in part by National Natural Science Foundation of China (No.: 12071174).

%%%%%%%%%%%%%%%%%%%%%%%%%%%%%%%%%%%%%%%%%%
%%%%%%%%%%%%%%%%%%%%%%%%%%%%%%%%%%%%%%%%%%
%%%%%%%%%%%%%%%%%%%%%%%%%%%%%%%%%%%%%%%%%%
%%%%%%%%%%%%%%%%%%%%%%%%%%%%%%%%%%%%%%%%%%
%%%%%%%%%%%%%%%%%%%%%%%%%%%%%%%%%%%%%%%%%%

%%%%%%%%%%%%%%%%%%%%%%%%%%%%%%%%%%%%%%%%%%

%%%%%%%%%%%%%%%%%%%%%%%%%%%%%%%%%%%%%%%%%%

%%%%%%%%%%%%%%%%%%%%%%%%%%%%%%%%%%%%%%%%%%

%%%%%%%%%%%%%%%%%%%%%%%%%%%%%%%%%%%%%%%%%%

%%%%%%%%%%%%%%%%%%%%%%%%%%%%%%%%%%%%%%%%%%

%%%%%%%%%%%%%%%%%%%%%%%%%%%%%%%%%%%%%%%%%%

%%%%%%%%%%%%%%%%%%%%%
% Bibliography
%%%%%%%%%%%%%%%%%%%%%

\vskip 1 cm

\bibliographystyle{plain}

\end{document}